\documentclass[12pt]{amsart}

\usepackage{amscd,amssymb,epsfig}

\def\Z{{\mathbb Z}}
\def\Q{{\mathbb Q}}
\def\R{{\mathbb R}}
\def\C{{\mathbb C}}
\def\P{{\mathbb P}}
\def\H{{\mathbb H}}
\def\V{{\mathbb V}}
\def\bA{{\mathbb A}}

\def\E{{\mathcal E}}
\def\M{{\mathcal M}}
\def\O{{\mathcal O}}
\def\X{{\mathcal X}}
\def\cC{{\mathcal C}}

\def\cL{{\mathcal L}}
\def\cP{{\mathcal P}}
\def\cR{{\mathcal R}}

\def\G{{\Gamma}}
\def\D{{\Delta}}
\def\S{{\Sigma}}
\def\L{{\Lambda}}

\def\a{{\alpha}}
\def\b{{\beta}}

\def\w{{\omega}}

\def\s{{\mathfrak s}}
\def\l{{\mathfrak l}}

\def\sl{{\s\l}}

\def\orb{{\mathrm{orb}}}

\def\Hbar{\overline{\H}}
\def\Mbar{\overline{\M}}
\def\Xbar{{\overline{X}}}
\def\Xtilde{\widetilde{X}}

\def\gtop{g_{\mathrm{top}}}
\def\Mtilde{\widetilde{M}}
\def\xtilde{\tilde{x}}
\def\piorb{\pi^{\mathrm{orb}}}

\def\dot{{\bullet}}
\def\bs{\backslash}

\def\blank{\phantom{x}}
\def\bil#1#2{\langle #1,#2 \rangle}
\def\intn{\langle \blank, \blank \rangle}

\renewcommand\Im{\operatorname{Im}}
\renewcommand\Re{\operatorname{Re}}
\newcommand\im{\operatorname{im}}               

\newcommand\tr{\operatorname{tr}}
\newcommand\Ad{\operatorname{Ad}}
\newcommand\rank{\operatorname{rank}}
\newcommand\Diff{\operatorname{Diff}}
\newcommand\Hom{\operatorname{Hom}}

\newcommand\Aut{\operatorname{Aut}}
\newcommand\Inn{\operatorname{Inn}}
\newcommand\Out{\operatorname{Out}}
\newcommand\Jac{\operatorname{Jac}}
\newcommand\Pic{\operatorname{Pic}}

\def\Picorb{{\Pic_{\mathrm{orb}}}}

\newtheorem{theorem}{Theorem}[section]

\newtheorem{proposition}[theorem]{Proposition}
\newtheorem{corollary}[theorem]{Corollary}

\theoremstyle{definition}
\newtheorem{definition}[theorem]{Definition}

\newtheorem{xca}[theorem]{Exercise}

\theoremstyle{remark}
\newtheorem{remark}[theorem]{Remark}

\begin{document}

\title[Moduli of Riemann Surfaces, Transcendental Aspects]{Moduli of Riemann
Surfaces, Transcendental Aspects}

\author{Richard Hain}

\address{Department of Mathematics, Duke University, Durham, NC 27708}

\email{hain@math.duke.edu}

\thanks{Please send comments and corrections to me at
{\tt hain@math.duke.edu}.}

\date{\today}

\maketitle


These notes are an informal introduction to moduli spaces of compact Riemann
surfaces via complex analysis, topology and Hodge Theory. The prerequisites
for the first lecture are just basic complex variables, basic Riemann surface
theory up to at least the Riemann-Roch formula, and some algebraic topology,
especially covering space theory. Some good references for this material
include \cite{ahlfors} for complex analysis, \cite{farkas-kra} and
\cite{forster} for the basic theory of Riemann surfaces, and \cite{greenberg}
for algebraic topology. For later lectures I will assume more. The book by
Clemens \cite{clemens} and Chapter~2 of Griffiths and Harris \cite{griff-harris}
are excellent and are highly recommended. Other useful references include the
surveys \cite{harer:survey} and \cite{hain-looijenga} and the book
\cite{harris-morrison}.

The first lecture covers moduli in genus 0 and genus 1 as these can be
understood using relatively elementary methods, but illustrate many of the
points which arise in higher genus. The notes cover more material than
was covered in the lectures, and sometimes the order of topics in the notes
differs from that in the lectures. I hope to add the material from the last
lecture on the Torelli group and Morita's approach to the tautological
classes in a future version.

\noindent{\it Acknowledgements.} 
It is a great pleasure to thank all those with whom I have had helpful
discussions while preparing these notes. I would especially like to thank
Robert Bryant, John Harer, Eduard Looijenga, Makoto Matsumoto and Mark
Stern.

\begin{center}
{\sc Lecture 1: Low Genus Examples}
\end{center}

Suppose that $g$ and $n$ are non-negative integers. An {\it $n$-pointed
Riemann surface $(C;x_1,\dots,x_n)$ of genus $g$} is a compact Riemann
surface $C$ of genus $g$ together with an ordered $n$-tuple of distinct
points $(x_1,\dots,x_n)$ of $C$. Two $n$-pointed Riemann surfaces
$(C;x_1,\dots,x_n)$ and $(C';x_1',\dots,x_n')$ are {\it isomorphic} if there
is a biholomorphism $f : C \to C'$ such that $f(x_j) = x_j'$ when
$1\le j \le n$.
The principal objects of study in these lectures are the spaces
$$
\M_{g,n} = \left\{
\parbox{2.75in}{isomorphism classes of $n$-pointed compact Riemann surfaces
$C$ of genus $g$} \right\}
$$
At the moment all we can say is that these are sets. One of the main
objectives of these lectures is to show that each $\M_{g,n}$ is a complex
analytic variety with very mild singularities.

Later we will only consider $\M_{g,n}$ when the stability condition
\begin{equation}
\label{stab_cond}
2g - 2 + n > 0
\end{equation}
is satisfied. But for the time being we will consider all possible values
of $g$ and $n$.
When $n=0$, we will simply write $\M_g$ instead of $\M_{g,0}$.

The space $\M_{g,n}$ is called the {\it moduli space of $n$-pointed curves
(or Riemann surfaces) of genus $g$}.
The isomorphism class of $(C;x_1,\dots,x_n)$ is called the
{\it moduli point} of $(C;x_1,\dots,x_n)$ and will be
denoted by $[C;x_1,\dots,x_n]$.

There are (at least) two notions of the genus of a compact Riemann
surface $C$. First there is the (analytic) genus
$$
g(C) := \dim H^0(C,\Omega^1_C),
$$
the dimension of the space of global holomorphic 1-forms on $C$. Second
there is the {\it topological genus}
$$
\gtop(C) := \frac{1}{2} \rank H_1(C,\Z).
$$
Intuitively, this is the `number of holes' in $C$. A basic fact
is that these are equal. There are various ways to prove this, but
perhaps the most standard is to use the Hodge Theorem (reference) which
implies that
$$
H^1(C,\C) \cong \{\text{holomorphic 1-forms}\} \oplus
\{\text{anti-holomorphic 1-forms}\}.
$$
The equality of $\gtop(C)$ and $g(C)$ follows immediately as complex
conjugation interchanges the holomorphic and antiholomorphic differentials.

Finally, we shall use the terms ``complex curve'' and ``Riemann surface''
interchangeably.

\section{Genus $0$}

It follows from Riemann-Roch formula that if $X$ is a compact Riemann
surface of genus $0$, then $X$ is biholomorphic to the Riemann sphere
$\P^1$. So $\M_0$ consists of a single point.

An {\it automorphism} of a Riemann surface $X$ is simply a biholomorphism
$f : X \to X$. The set of all automorphisms of $X$ forms a group $\Aut X$.
The group $GL_2(\C)$ acts in $\P^1$ via fractional linear transformations:
$$
\begin{pmatrix}
a & b \cr c & d
\end{pmatrix}
: z \mapsto \frac{az+b}{cz+d}
$$
The scalar matrices $S$ act trivially, and so we have a homomorphism
$$
PGL_2(\C) \to \Aut \P^1
$$
where for any field $F$
$$
PGL_n(F) = GL_n(F)/\text{\{scalar matrices\}}
$$
and
$$
PSL_n(F) = SL_n(F)/\text{\{scalar matrices of determinant 1\}}.
$$

\begin{xca}
Prove that $PGL_2(\C) \cong PSL_2(\C)$ and that these are isomorphic
to $\Aut \P^1$.
\end{xca}

\begin{xca}
Prove that $\Aut \P^1$ acts 3-transitively on $\P^1$. That is, given
any two ordered 3-tuples $(a_1,a_2,a_3)$ and $(b_1,b_2,b_3)$ of distinct 
points of $\P^1$, there is an element $f$ of $\Aut\P^1$ such that
$f(a_j) = b_j$ for $j=1,2,3$. Show that $f$ is unique.
\end{xca}

\begin{xca}
Prove that if $X$ is a compact Riemann surface of genus $0$, then $X$
is biholomorphic to the Riemann sphere.
\end{xca}

\begin{xca}
Show that the automorphism group of an $n$-pointed curve of genus $g$
is finite if and only if the stability condition~(\ref{stab_cond}) is
satisfied. (Depending on what you know, you may find this a little 
difficult at present. More techniques will become available soon.)
\end{xca}

Since $\Aut \P^1$ acts 3-transitively on $\P^1$, we have:

\begin{proposition}
Every $n$-pointed Riemann surface of genus 0 is isomorphic to
$$
\begin{matrix}
(\P^1;\infty) & \text{ if } n=1; \cr
(\P^1;0,\infty) & \text{ if } n=2; \cr
(\P^1;0,1,\infty) & \text{ if } n= 3.
\end{matrix}
$$ \qed
\end{proposition}

\begin{corollary}
If $0 \le n\le 3$, then $\M_{0,n}$ consists of a single point. \qed
\end{corollary}

The first interesting case is when $n=4$. If $(X;x_1,x_2,x_3,x_4)$ is a
4-pointed Riemann surface of genus 0, then there is a unique biholomorphism
$f : X \to \P^1$ with $f(x_2) = 1$, $f(x_3) = 0$ and $f(x_4) = \infty$.
The value of $f(x_1)$ is forced by these conditions. Since the $x_j$ are
distinct and $f$ is a biholomorphism, $f(x_1) \in \C - \{0,1\}$. It is 
therefore an invariant of $(X;x_1,x_2,x_3,x_4)$.

\begin{xca}
Show that if $g : X \to \P^1$ is any biholomorphism, then $f(x_1)$ is the
cross ratio
$$
(g(x_1):g(x_2):g(x_3):g(x_4))
$$
of $g(x_1)$, $g(x_2)$, $g(x_3)$, $g(x_4)$. Recall that the cross ratio of 
four distinct points $x_1$, $x_2$, $x_3$, $x_4$ in $\P^1$ is defined by
$$
(x_1:x_2:x_3:x_4) = \frac{(x_1-x_3)/(x_2-x_3)}{(x_1-x_4)/(x_2-x_4)}
$$
\end{xca}

The result of the previous exercise can be rephrased as a statement about
moduli spaces:

\begin{proposition}
The moduli space $\M_{0,4}$ can be identified naturally with $\C - \{0,1\}$.
The moduli point $[\P^1;x_1,x_2,x_3,x_4]$ is identified with the cross
ratio $(x_1:x_2:x_3:x_4) \in \C-\{0,1\}$. \qed
\end{proposition}

It is now easy to generalize this to general $n \ge 4$. Since every genus
0 Riemann surface is biholomorphic to $\P^1$, we need only consider $n$-pointed
curves of the form $(\P^1;x_1,\dots,x_n)$. There is a unique automorphism
$f$ of $\P^1$ such that $f(x_1) = 0, f(x_2) = 1$ and $f(x_3) = \infty$.
So every $n$-pointed Riemann surface of genus 0 is isomorphic to exactly
one of the form
$$
(\P^1;0,1,\infty,y_1,\dots,y_{n-3}).
$$
To say that this is an $n$-pointed curve is to say that the points
$0, 1, \infty$, $y_1,\dots,y_{n-3}$ are distinct. That is, it
$$
(y_1,\dots,y_{n-3}) \in (\C - \{0,1\})^{n-3} - \D
$$
where $\D = \cup_{j<k} \D_{jk}$ is the union of the diagonals
$$
\D_{jk} = \{(y_1,\dots,y_{n-3}): y_j = y_k\}.
$$
This is an affine algebraic variety as it is the complement of a
divisor in an affine space. This shows that:

\begin{theorem}
If $n\ge 3$, then $\M_{0,n}$ is a smooth affine algebraic variety of
dimension $n-3$ isomorphic to
$$
(\C - \{0,1\})^{n-3} - \D.
$$
\end{theorem}

The symmetric group $\S_n$ acts on $\M_{0,n}$ by
$$
\sigma :
(\P^1;x_1,\dots,x_n) \mapsto (\P^1;x_{\sigma(1)},\dots, x_{\sigma(n)}).
$$
 
\begin{xca}
Show that each $\sigma \in \S_n$ acts on $\M_{0,n}$ as a regular mapping.
(Hint: it suffices to consider the case of a transposition.)
\end{xca}

\begin{xca}
Suppose that $n\ge 3$. Construct a universal $n$-pointed genus 0 curve
$\M_{0,n}\times\P^1 \to \M_{0,n}$ that is equipped with $n$ disjoint sections
$\sigma_1,\dots,\sigma_n$ such that
$$
\sigma_j([\P^1;x_1,\dots,x_n]) = ([\P^1;x_1,\dots,x_n],x_j).
$$
Show that it is universal in the sense that if $f:X \to T$ is a family of
smooth genus 0 curves over a smooth variety $T$ and if the family has
$n$ sections $s_1,\dots,s_n$ that are disjoint, then there is a holomorphic
mapping $\phi_f : T \to \M_{0,n}$ such that the pullback of the universal
family is $f$ and the pullback of $\sigma_j$ is $s_j$.
\end{xca}


\section{Genus $1$}

The study of the moduli space of genus 1 compact Riemann surfaces is very
rich and has a long history because of its fundamental connections to number
theory and the theory of plane cubic curves. We will take a transcendental
approach to understanding $\M_1$ which will reveal the connection with
modular forms. Our first task is show that genus 1 Riemann surfaces can
always be represented as the quotient of $\C$ by a lattice.

One way to construct a Riemann surface of genus 1 is to take the quotient of
$\C$ by a lattice. Recall that a {\it lattice} in a finite dimensional real
vector space $V$ is a finitely generated (and therefore free abelian) subgroup
$\L$ of $V$ with the property that a basis of $\L$ as an abelian group is
also a basis of $V$ as a real vector space. A lattice in $\C$ is thus a
subgroup $\L$ of $\C$ that is isomorphic to $\Z^2$ and is generated by
two complex numbers that are not real multiples of each other.

\begin{xca}
Show that if $\L$ is a lattice in $V$ and if $\dim_\R V = d$, then
$V/\L$ is a compact manifold of real dimension $d$, which is diffeomorphic
to the $d$-torus $(\R/\Z)^d$.
\end{xca}

If $\L$ is a lattice in $\C$, the quotient group $\C/\L$ is a compact
Riemann surface which is diffeomorphic to the product of two circles, and
so of genus 1.

\begin{theorem}
\label{lattice}
If $C$ is a compact Riemann surface of genus $1$, then there
is a lattice $\L$ in $\C$ and an isomorphism $\mu : C \to \C/\L$.
If $x_o \in C$, then we may choose $\mu$ such that $\mu(x_o)=0$.
\end{theorem}

The proof follows from the sequence of exercises below. Let
$C$ be a compact Riemann surface of genus $1$.

\begin{xca}
Show that a non-zero holomorphic differential on $C$ has no zeros.
Hint: Use Riemann-Roch.
\end{xca}

Since $\gtop(C) = g(C) = 1$, we know that $H_1(C,\Z)$ is free of rank 2.
Fix a non-zero holomorphic differential $w$ on $C$. Every other holomorphic
differential is a multiple of $w$. The {\it period lattice} of $C$ is 
defined to be
$$
\L = \big\{\int_c w : c \in H_1(C,\Z)\big\}.
$$ 
This is easily seen to be a subgroup of $\C$.

\begin{xca}
Show that $\L$ is a lattice in $\C$. Hint: Choose a basis $a,b$ of $H_1(C,\Z)$.
Show that if $\int_a w$ and $\int_b w$ are linearly independent over $\R$,
then this would contradict the Hodge decomposition of $H^1(C,\C)$.
\end{xca}

Let $E = \C/\L$. Our next task is to construct a holomorphic mapping
from $C$ to $E$.

\begin{xca}
Fix a base point $x_o$ of $C$. Define a mapping $\nu : C \to E$ by
$$
\nu(x) = \int_\gamma w
$$
where $\gamma$ is any smooth path in $C$ that goes from $x_o$ to $x$.
Show that
\begin{enumerate}
\item $\nu$ is well defined;
\item $\nu$ is holomorphic;
\item $\nu$ has nowhere vanishing differential, and is therefore a covering
map;
\item the homomorphism $\nu_\ast : \pi_1(C,x_o) \to \pi_1(E,0) \cong \Lambda$
is surjective, and therefore an isomorphism.
\end{enumerate}
Deduce that $\nu$ is a biholomorphism.
\end{xca}

This completes the proof of Theorem~\ref{lattice}. It has the following
important consequence:

\begin{corollary}
\label{transl}
If $C$ is a compact Riemann surface of genus $1$, then the automorphism
group of $C$ acts transitively on $C$. Consequently, the natural mapping
$\M_{1,1} \to \M_1$ that takes $[C;x]$ to $[C]$ is a bijection.
\end{corollary}

\begin{proof}
This follows as every genus 1 Riemann surface is isomorphic to one
of the form $\C/\L$. For such Riemann surfaces, we have the homomorphism
$$
\C/\L \to \Aut (\C/\L)
$$
that takes the coset $a + \L$ to the translation $z + \L \mapsto z+a+\L$.
\end{proof}

\begin{definition}
An {\it elliptic curve} is a a genus 1 curve $C$ together with a point
$x_o\in C$.
\end{definition}

The previous result says that if $C$ is a genus 1 curve and $x_o$ and $y_o$
are points of $C$, then the elliptic curves $(C;x_o)$ and $(C;y_o)$ are
isomorphic.

The moduli space of elliptic curves is $\M_{1,1}$.

\begin{xca}
Suppose that $f:C \to \C/\L$ is a holomorphic mapping from an arbitrary
Riemann surface to $\C/\L$. Let $x_o$ be a base point of $C$.
The 1-form $dz$ on $\C$ descends to a holomorphic differential $w$ on
$\C/\L$. Its pullback $f^\ast w$ is a holomorphic differential on $C$.
Show that for all $x\in C$,
$$
f(x) = f(x_o) + \int_\gamma f^\ast w + \L
$$
where $\gamma$ is a path in $C$ from $x_o$ to $x$.
\end{xca}

\begin{xca}
Use the results of the previous exercise to prove the following result.
\end{xca}

\begin{corollary}
\label{rescale}
If $\L_1$ and $\L_2$ are lattices in $\C$, then
$\C/\L_1$ is isomorphic to $\C/\L_2$ if and only there exists
$\lambda \in \C^\ast$ such that $\L_1 = \lambda \L_2$. \qed
\end{corollary}

\begin{xca}
Show that if $(C;x_o)$ is an elliptic curve, then $C$ has a natural group
structure with identity $x_o$.
\end{xca}

We are finally ready to give a construction of $\M_1$.
Recall that the complex structure on a Riemann surface $C$ gives it a 
canonical orientation. This can be thought of as giving a direction of
``positive rotation'' about each point in the surface --- the
positive direction being that given by turning counter-clockwise about
the point in any local holomorphic coordinate system. If $C$ is compact,
this orientation allows us to define the intersection number of two
transversally intersecting closed curves on $C$. It depends only on the
homology classes of the two curves and therefore defines the
{\it intersection pairing}
$$
\intn : H_1(C,\Z) \otimes H_1(C,\Z) \to \Z.
$$
If $\a$, $\b$ is a basis of $H_1(C,\Z)$, then $\bil \a \b = \pm 1$. We shall
call the basis {\it positive} if $\bil \a \b = 1$.

A {\it framing} of Riemann surface $C$ of genus 1 is a positive basis $\a,\b$
of its first homology group. We will refer to $(C:\a,\b)$ as a {\it framed}
genus 1 Riemann surface. Two framed genus 1 Riemann surfaces $(C:\a,\b)$ and
$(C':\a',\b')$ are isomorphic if there is a biholomorphism $f:C\to C'$ such
that $\a' = f_\ast\a$ and $\b' = f_\ast \b$.

Let
$$
\X_1 = \left\{
\parbox{2.25in}{isomorphism classes of framed Riemann surfaces of genus 1
} \right\}.
$$
At the moment, this is just a set. But soon we will see that it is itself
a Riemann surface. Note that forgetting the framing defines a
a function
$$
\phi : \X_1 \to \M_1.
$$
Denote the isomorphism class of $(C:\a,\b)$ by $[C:\a,\b]$.
If $C$ is a genus 1 Riemann surface, then
$$
\phi^{-1}([C]) =
\{[C:\a,\b] : (\a,\b) \text{ is a positive basis of }H_1(C,\Z)\}.
$$

\begin{xca}
Show that if $(\a,\b)$ and $(\a',\b')$ are two positive bases of $H_1(C,\Z)$,
then there is a unique element
$$
\begin{pmatrix} a & b \cr c & d \end{pmatrix}
$$
of $SL_2(\Z)$ such that
$$
\begin{pmatrix} \b' \cr \a' \end{pmatrix} =
\begin{pmatrix} a & b \cr c & d \end{pmatrix}
\begin{pmatrix} \b \cr \a \end{pmatrix}
$$
(The reason for writing the basis vectors in the reverse order will become
apparent shortly.)
\end{xca}

Define an action of $SL_2(\Z)$ on $\X_1$ by
$$
\begin{pmatrix} a & b \cr c & d \end{pmatrix}[C:\a,\b] = [C:\a',b']
$$
where
$$
\begin{pmatrix} \b' \cr \a' \end{pmatrix} =
\begin{pmatrix} a & b \cr c & d \end{pmatrix}
\begin{pmatrix} \b \cr \a \end{pmatrix}
$$

\begin{xca}
\label{quot}
Show that there is a natural bijection
$$
\M_1 \cong SL_2(\Z)\bs\X_1.
$$
\end{xca}

At present, $\X_1$ is just a set, but we now show that it is naturally
a Riemann surface. We know from Theorem~\ref{lattice} that every element
of $\X_1$ is of the form $[\C/\L:\a,\b]$. But, by standard algebraic
topology, there is a natural isomorphism
$$
\L \cong H_1(C,\Z).
$$
Thus a basis of $H_1(C,\Z)$ corresponds to a basis of $\L$.

\begin{xca}
Suppose that $\a,\b$ is a basis of $H_1(\C/\L,\Z)$ and that $\w_1,\w_2$
is the corresponding basis of $\L$. Show that $\a,\b$ is positive if and
only if $\w_2/\w_1$ has positive imaginary part.
\end{xca}

It follows from this and Corollary~\ref{rescale} that
$$
\X_1 = \{[\C/\L:\w_1,\w_2]: \Im (\w_2/\w_1) > 0 \}/\C^\ast
$$
where the $\C^\ast$-action is defined by
$$
\lambda\cdot [\C/\L:\w_1,\w_2] = [\C/\lambda\L:\lambda\w_1,\lambda\w_2].
$$
We can go even further: since the basis $\w_1$, $\w_2$ determines the
lattice,
$$
\L = \Z\w_1 \oplus \Z\w_2,
$$
we can dispense with the lattice altogether. We have:
\begin{equation}
\label{last}
\X_1 =
\{(\w_1,\w_2): \w_1,\w_2 \in \C \text{ and }\Im (\w_2/\w_1) > 0 \}/\C^\ast
\end{equation}
where $\C^\ast$ acts on $(\w_1,\w_2)$ by scalar multiplication.

Denote the upper half plane $\{z\in \C : \Im z > 0\}$ by $\H$. Each
$\tau \in \H$ determines the element
$$
[\C/(\Z \oplus \Z\tau):1,\tau]
$$
of $\X_1$. This defines a function $\psi : \H \to \X_1$. Under the
identification (\ref{last}),
$$
\psi(\tau) = \text{ the $\C^\ast$-orbit of }(1,\tau)
$$
Since $(\w_1,\w_2)$ and $(1,\w_2/\w_1)$ are in the same orbit, we have
proved:

\begin{theorem}
\label{bijn}
The function $\psi : \H \to \X_1$ is a bijection. \qed
\end{theorem}

The group $PSL_2(\C)$ acts on $\P^1$ by fractional linear transformations.

\begin{xca}
Show that $T \in PSL_2(\C)$ satisfies $T(\H) \subseteq \H$ if and only if
$T \in PSL_2(\R)$.
\end{xca}

Thus the group $SL_2(\R)$ acts on $\H$ by fractional linear transformations:
$$
\begin{pmatrix} a & b \cr c & d \end{pmatrix} \tau =
\frac{a\tau + b}{c\tau + d}
$$

\begin{xca}
\label{equivariant}
Show that $\psi : \H \to \X_1$ is $SL_2(\Z)$-equivariant. That is, if
$T \in SL_2(\Z)$, then $T\psi(\tau) = \psi(T\tau)$.
\end{xca}

\begin{theorem}
\label{M1}
There are natural bijections
$$
\M_1 \cong \M_{1,1} \cong SL_2(\Z)\bs \H.
$$
\end{theorem}

\begin{proof}
The first bijection was established in Corollary~\ref{transl}. The second
follows from Exercise~\ref{quot}, Theorem~\ref{bijn} and
Exercise~\ref{equivariant}.
\end{proof}

\begin{xca}
Suppose that $C$ is a genus 1 Riemann surface. Show that the point
of $SL_2(\Z)\bs \H$ that corresponds to $[C]\in \M_1$ is the $SL_2(\Z)$
orbit of
$$
\bigg(\int_\b w \bigg/ \int_\a w\bigg) \in \H
$$ 
where $w$ is any non-zero element of $H^0(C,\Omega^1_C)$ and $\a$, $\b$
is a positive basis of $H_1(C,\Z)$.
\end{xca}

\subsection{Understanding $SL_2(\Z)\bs \H$}

The Riemann surface structure on $\H$ descends to a Riemann surface
structure on $SL_2(\Z)\bs \H$. Good references for this are Chapter~VII
of Serre's book \cite{serre}, and Chapter~3 of Clemens' book \cite{clemens}.
We'll sketch part of the proof of the following fundamental theorem.

\begin{theorem}
The quotient of $\H$ by $SL_2(\Z)$ has a unique structure of a Riemann
surface such that the projection $\H \to SL_2(\Z) \bs \H$ is holomorphic.
Moreover, there is a biholomorphism between $SL_2(\Z) \bs \H$ and $\C$
which can be given by the modular function $j:\H \to \C$, where
$$
j(\tau) = \frac{1}{q} + 744 + 196\,884\,q + 21\,493\,760\, q^2 + \cdots
$$
and $q = e^{2\pi i \tau}$.
\end{theorem}

The following exercises will allow you to construct most of the
proof. The rest can be found in \cite{serre} and \cite{clemens}.

Let $\P^1(\R) = \R \cup \{\infty\}$. This is a circle on the Riemann
sphere which forms the boundary of $\H$. Let $\Hbar$ be the closure of
$\H$ in the Riemann sphere $\P^1$; it is the union of $\H$ and $\P^1(\R)$.
Recall that every non-trivial element of $PSL_2(\C)$ has at most two
fixed points in $\P^1$. Note that the fixed points of elements of $PSL_2(\R)$
are real or occur in complex conjugate pairs.

\begin{xca}
Suppose that $T \in SL_2(\Z)$ is not a scalar matrix. Show that $T$ has
exactly
\begin{enumerate}
\item one fixed point in $\H$ if and only if $|\tr T| < 2$;
\item one fixed point in $\P^1(\R)$ if and only if $|\tr T| = 2$;
\item two fixed points in $\P^1(\R)$ if and only if $|\tr T| > 2$.
\end{enumerate}
Show that $T\in SL_2(\R)$ has finite order if and only if $T$ has a fixed
point in $\H$.
\end{xca}

Fix an integer $l\ge 0$. The {\it level $l$ subgroup} of $SL_2(\Z)$ is
the subgroup of $SL_2(\Z)$ consisting of those matrices congruent to
the identity mod $l$. We shall denote it by $SL_2(\Z)[l]$. Since it is
the kernel of the homomorphism $SL_2(\Z) \to SL_2(\Z/l)$, it is normal
and of finite index in $SL_2(\Z)$.

\begin{xca}
Show that $SL_2(\Z)[l]$ is torsion free for all $l\ge 3$. Hint: use the
previous exercise. Deduce that $SL_2(\Z)[l]\bs \H$ is a Riemann surface
with fundamental group $SL_2(\Z)[l]$ and universal covering $\H$ 
whenever $l\ge 3$.
\end{xca}

The quotient of $SL_2(\Z)$ by its level $l$ subgroup is $SL_2(\Z/l)$,
which is a finite group. It follows that the projection
$$
SL_2(\Z)[l]\bs \H \to SL_2(\Z)\bs \H
$$
is finite-to-one of degree equal to the half order of $SL_2(\Z/l)$
when $l>2$.

\begin{xca}
Show that the quotient $G\bs C$ of a Riemann surface $C$ by a finite subgroup
$G$ of $\Aut C$ has the structure of a Riemann surface such that the
projection
$C \to G \bs C$ is holomorphic. Hint: first show that for each point $x$ of
$C$, the isotropy group
$$
G_x := \{g \in G : gx=x\}
$$
is cyclic. This can be done by considering the actions of $G_x$ on $T_x X$
and $\O_{X,x}$.
\end{xca}

This result, combined with the previous exercises, establishes that
$\M_1 = SL_2(\Z)\bs \H$ has a natural structure of a Riemann surface
such that the projection $\H \to \M_1$ is holomorphic.

Serre \cite[p.~78]{serre} proves that a fundamental domain for the action
of $SL_2(\Z)$ on $\H$ is the region
$$
F = \{\tau \in \H : |\Re \tau| \le 1/2, |\tau| \ge 1 \}.
$$
\begin{figure}[!ht]
\epsfig{file=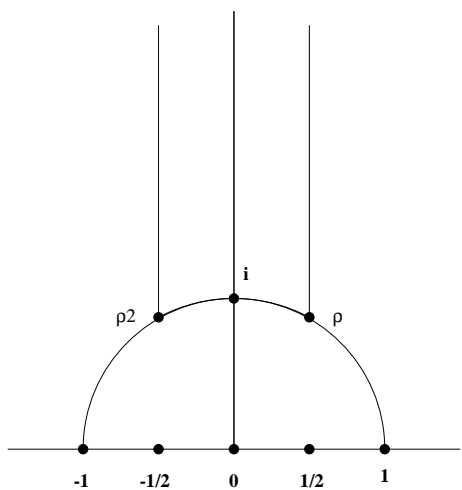, height=2.5in}
\caption{A
a fundamental domain of $SL_2(\Z)$ in $\H$}
\end{figure}
Points of $F$ can be thought of as giving a canonical framing
of a lattice $\Lambda$ in $\C$. Such a framing is given as follows:
the first basis element is a non-zero vector of shortest length in
$\Lambda$, the second basis element is a shortest vector in $\Lambda$
that is not a multiple of the first.

\begin{xca}
Show that $\tau \in F$ if and only if $(1,\tau)$ is such a canonical
basis of the lattice $\Z \oplus \Z\tau$.
\end{xca}

Serre \cite[p.~78]{serre} also proves that $PSL_2(\Z)$ is generated by
$\tau \mapsto -1/\tau$ and $\tau \mapsto \tau + 1$.

\begin{xca}
Use this to prove that $SL_2(\Z)\bs \H$ is the quotient of the fundamental
domain $F$ obtained by identifying the opposite vertical sides, and by
identifying the arc of the circle $|\tau|=1$ from $\rho$ to $i$ with
the arc from $\rho^2$ to $i$. Deduce that $\M_1=SL_2(\Z)\bs \H$ is a
Riemann surface homeomorphic to a disk.
\end{xca}

It is conceivable that $\M_1$ is biholomorphic to a disk, for example. But
this is not the case as $\M_1$ can be compactified by adding one point.

\begin{xca}
Show that $\M_1$ can be compactified by adding a single point
$\infty$. A coordinate neighbourhood of $\infty$ is the unit disk $\D$.
Denote the holomorphic coordinate in it by $q$. The point $\tau$ of
$\M_1 = SL_2(\Z)\bs \H$ is identified with the point $e^{2\pi i\tau}$ of
$\D$. Show that
$$
\M_1 \cup \{\infty\}
$$
is a compact Riemann surface of genus 0 where $q$ is a local parameter
about $\infty$ and where $\M_1$ holomorphically embedded. Deduce that
$\M_1$ is biholomorphic to $\C$.
\end{xca}

\begin{remark}
This compactification is the moduli space $\Mbar_{1,1}$.
\end{remark}

Since every compact Riemann surface is canonically a complex algebraic
curve, this shows that $\M_1$ is an algebraic variety.

\subsection{Automorphisms}

The automorphisms of an elliptic curve are intimately related with the
set of elements of $SL_2(\Z)$ that stabilize the points corresponding
to it in $\H$.

\begin{xca}
Suppose that $C$ is a genus 1 Riemann surface curve. Suppose that
$\phi : C \to C$ is an automorphism of $C$. Show that 
\begin{enumerate}
\item $\phi$ is a translation if and only if
$\phi_\ast : H_1(C,\Z) \to H_1(C,\Z)$ is the identity;
\item if $\a,\b$ is a positive basis of $H_1(C,\Z)$, then
$$
[C:\a,\b] = [C: \phi_\ast \a, \phi_\ast \b] \in \X_1;
$$
\end{enumerate}
Deduce that there is a natural isomorphism
$$
\Aut (C;x_o) \cong \text{ stabilizer in $SL_2(\Z)$ of }[C:\a,\b] \in \X_1.
$$
\end{xca}

Note that any two points of $\H$ that lie in the same orbit of $SL_2(\Z)$
have isomorphic stabilizers. Consequently, to find all elliptic curves with
automorphism groups larger than $\Z/2\Z$, one only has to look for points
in the fundamental domain with stabilizers larger than $\Z/2\Z$.

\begin{xca}
Show that the stabilizer in $SL_2(\Z)$ of $\tau \in \H$ is
$$
\begin{cases}
\Z/2\Z & \tau \notin \text{ orbit of $i$ and $\rho$}; \cr
\Z/4\Z & \tau \in \text{ orbit of $i$}; \cr
\Z/6\Z & \tau \in \text{ orbit of $\rho$}.
\end{cases}
$$
\end{xca}

An immediate consequence of this computation is the following:

\begin{theorem}
If $(C;x_o)$ is an elliptic curve, then
$$
\Aut(C,x_o) \cong
\begin{cases}
\Z/4\Z & \text{ if } (C;x_o) \cong (\C/\Z[i];0) \cr
\Z/6\Z & \text{ if } (C;x_o) \cong (\C/\Z[\rho];0) \cr
\Z/2\Z & \text{ otherwise.}
\end{cases}
$$ 
\end{theorem}

It is easy to see from this, for example, that the Fermat cubic
$$
x^3 + y^3 + z^3 = 0
$$
in $\P^2$ has automorphism group $\Z/6\Z$ and therefore is isomorphic
to $\Z/\Z[\rho]$.

\subsection{Families of Genus 1 and Elliptic Curves}

Suppose that $X$ is a complex analytic manifold and that $f:X \to T$ is a
holomorphic mapping to another complex manifold each of whose
fibers is a genus 1 curve. We say that $f$ is a family of genus 1 curves. If
there is a section $s:T \to X$ of $f$, then the fiber of $f$ over $t\in T$ is
an elliptic curve with identity $s(t)$. We say that $f$ is a family of
elliptic curves.

Each such family gives rise to a function
$$
\phi_f : T \to \M_1
$$
that is defined by taking $t\in T$ to the moduli point $[X_t]$ of the
fiber $X_t$ of $f$ over $t$. We shall call it the {\it period mapping}
of the family.

\begin{theorem}
The period map $\phi_f$ is holomorphic.
\end{theorem}

\begin{proof}[Sketch of Proof]
For simplicity, we suppose that $T$ has complex dimension 1. The relative
holomorphic tangent bundle of $f$ is the holomorphic line bundle
$T'_f$ on $X$ consisting of holomorphic tangent vectors to $X$ that
are tangent to the fibers of $f$. That is,
$$
T_f' = \ker\{T'X \to T'T\}.
$$
The sheaf of holomorphic sections of its dual is called the relative
dualizing sheaf and is often denoted by $\w_{X/T}$. The push forward
$f_\ast \w_{X/T}$ of this sheaf to $T$ is a holomorphic line bundle
over $T$ and has fiber
$$
H^0(X_t,\Omega^1_{X_t})
$$
over $t\in T$. Fix a reference point $t_o \in T$. Let $w(t)$ be a local
holomorphic section of $\w_{X/T}$ defined in a contractible neighbourhood
$U$ of $t_o$. Since bundle is locally topologically trivial,  $f^{-1}(U)$
is homeomorphic to $U \times X_{t_o}$. We can thus identify $H_1(X_t,\Z)$
with $H_1(X_{t_o},\Z)$ for each $t\in U$. Fix a positive basis $\a,\b$ of
$H_1(X_{t_o},\Z)$. This can be viewed as a positive basis of $H_1(X_t,\Z)$
for all $t\in U$.

It is not difficult to show that
$$
\int_\a w(t) \text{ and } \int_\b w(t)
$$
vary holomorphically with $t\in U$. It follows that
$$
\tau(t) := \bigg\{\int_\b w(t) \bigg/ \int_\a w(t) \bigg\}
$$
varies holomorphically with $t\in U$. This shows that map $\phi_f$ is
holomorphic in the neighbourhood $U$ of $t_o$. It follows that $\phi_f$
is holomorphic.
\end{proof}

If you examine the proof, you will see that we really proved two extra facts.
First, the period mapping $\phi_f : T \to \M_1$ associated to every
family $f : X \to T$ of genus 1 curves is {\it locally liftable} to a
holomorphic mapping to $\H$. If, in addition, $f$ is a family of elliptic
curves, then the period mapping $\phi_f$ determines the family $f$ and
the section $s$:

\begin{proposition}
\label{class_map}
The period mapping $\phi_f$ associated to a family $f : X \to T$ of genus $1$
curves is locally liftable to a holomorphic mapping to $\H$. If $f$ is a
family of elliptic curves, then $\phi_f$ can be globally lifted to a holomorphic
mapping $\tilde{\phi}_f : \tilde{T} \to \H$ and there is a homomorphism
$\phi_{f\ast} : \pi_1(T,t_o) \to SL_2(\Z)$ (unique up to conjugacy) such that
the diagram
$$
\begin{CD}
\tilde{T} @>\tilde{\phi}_f>> \H \cr
@VVV @VVV \cr
T @>>\phi_f> \M_1
\end{CD}
$$
commutes and such that
$$
\tilde{\phi}_f(\gamma\cdot x) = \phi_{f\ast}(\gamma)\cdot \tilde{\phi}_f(x)
$$
for all $x\in \tilde{T}$ and all $\gamma \in \pi_1(T,t_o)$. \qed
\end{proposition}

The proof is left as an exercise. We shall see in a moment that the 
converse of this result is also true.

This result has an important consequence --- and that is that not every
holomorphic mapping $T \to \M_1$ is the period mapping of a holomorphic
family of elliptic (or even genus 1) curves. The reason for this is that not
every mapping $T \to \M_1$ is locally liftable.

\begin{xca}
Show that the identity mapping $\M_1 \to \M_1$ is not locally liftable.
Deduce that there is no family of genus 1 curves over $\M_1$ whose period
mapping is the identity mapping. In particular, show that there is no universal
elliptic curve over $\M_1$.
\end{xca}

To each genus 1 curve $C$, we can canonically associate the elliptic curve
$\Jac C := \Pic^0 C$ which we shall call the {\it jacobian} of $C$. Abel's
Theorem tells us that $C$ and $\Jac C$ are isomorphic as genus 1 curves, but
the isomorphism depends on the choice of a base point of $C$.

\begin{xca}
Show that for each family $f:X \to T$ of genus 1 curves the corresponding
family of jacobians
is a family of elliptic curves. Show that these two families have the
same period mapping. Show that if $f$ is a family of elliptic curves, then
the family of jacobians is canonically isomorphic to the original family $f$.
\end{xca}

\section{Orbifolds}

The discussion in the previous section suggests that $\M_1$ should be
viewed as $SL_2(\Z) \bs \H$ rather than as $\C$ --- or that $\M_1$ is
not an algebraic variety or a manifold, but rather something whose local
structure includes the information of how it is locally the quotient of a
disk by a finite group. In topology such objects are called {\it orbifolds}
and in algebraic geometry {\it stacks}. Very roughly speaking, orbifolds are
to manifolds as stacks are to varieties. The moduli spaces $\M_{g,n}$ are
often conveniently viewed as orbifolds or as stacks.

For us, an orbifold is a topological space that is the quotient of a simply
connected topological space $X$ by a group $\G$. The group is required to
act properly discontinuously on $X$ and all isotropy groups are required to be
finite.\footnote{We could have added some other natural
conditions, such as requiring that there be
a finite index subgroup of $\G$ that acts on $X$ fixed point freely. We
could also require that the set of elements of $G$ that act trivially on
$X$ is central in $\G$. Both of these conditions are natural and are
satisfied by all of our primary examples of orbifolds, the $\M_{g,n}$.}
For example, $\M_1$ can be viewed as an orbifold as the quotient of $\H$ by
$SL_2(\Z)$.

This definition is not as general as it could be, but since all of our
moduli spaces
$\M_{g,n}$ are of this form, it is good enough for our purposes. The 
general definition is obtained by ``sheafifying'' this one --- i.e., 
general orbifolds are locally the quotient of a simply connected space
by a finite group. A general definition along these lines can be found in
Chapter~13 of \cite{thurston}. Mumford's definition of stacks can be found
in \cite{mumford}.

A morphism $f : \G_1 \bs X_1 \to \G_2 \bs X_2$ of orbifolds is a continuous
map that arises as follows: there is a homomorphism $f_\ast : \G_1 \to \G_2$
and a continuous mapping $\tilde{f} : X_1 \to X_2$ that is equivariant
with respect to $f_\ast$; that is, the diagram
$$
\begin{CD}
X_1 @>\tilde{f}>> X_2 \cr
@VgVV @VV{f_\ast(g)}V \cr
X_1 @>\tilde{f}>> X_2
\end{CD}
$$
commutes for all $g \in \G_1$.

At this stage, I should point out that every reasonable topological space
$X$ can be viewed as the quotient $\pi_1(X,x)\bs \tilde{X}$ of its universal
covering by its fundamental group. In this way, every topological space can
be regarded as an orbifold. It thus it makes sense to talk about orbifold
mappings between orbifolds and ordinary topological spaces.

An orbifold $\G\bs X$
can have an enriched structure --- such as a smooth, riemannian, K\"ahler,
or algebraic structure. One just insists that its  ``orbifold universal
covering'' $X$ has such a structure and that $\G$ act on $X$ as automorphisms
of this structure. Maps between two orbifolds with the same kind of enriched
structure are defined in the obvious way. For example, a map between two
orbifolds with complex structures is given by an equivariant holomorphic map
between their orbifold universal coverings.

Many orbifolds are given as quotients of non-simply connected spaces
by a group that acts discontinuously, but not fixed point freely. Such
quotients have canonical orbifold structures: If $M$ is a topological space
and $G$ a group that acts discontinuously on $M$, then
$$
G\bs M \cong \G \bs \Mtilde
$$
where $p:\Mtilde \to M$ is the universal covering\footnote{Of course, I am
assuming that $M$ is nice enough as a topological space to have a universal
covering. We assume, for example, that $M$ is locally simply connected, a
condition satisfied by all complex algebraic varieties and all manifolds.}
of $M$ and
$$
\G =
\{(\phi,g) : \text{ where }\phi :\Mtilde \to \Mtilde \text{ covers } g\in G\}.
$$
Here $\phi$ covers $g\in G$ means that the diagram
$$
\begin{CD}
\Mtilde @>\phi>> \Mtilde \cr
@VpVV @VVpV \cr
M @>>g> M
\end{CD}
$$
commutes.

There is a notion of the {\it orbifold fundamental} group
$\piorb_1(\G\bs X,x_o)$ of a connected pointed orbifold. It is a variant of
the definition of the fundamental group of a pointed topological space.

First, we'll fix the convention that the composition
$\a\b$ of two paths $\a$ and $\b$ in a space is defined when $\a(1)=\b(0)$.
The product path is the one obtained by traversing $\a$ first, then $\b$.
This is the convention used in \cite{greenberg}, for example, and
is the opposite of the convention used by many algebraic geometers such
as Deligne.

We need to impose several conditions on $X$ and on $x_o$ in order for the
definition to make sense.
Let $p : X \to \G\bs X$ be the projection. Note that if $x,y \in p^{-1}(x_o)$,
then the isotropy groups $\G_y$ and $\G_z$ are conjugate in $\G$. Our first
assumption is that $\G_y$ is contained in the center of $\G$ for one (and
hence all) $y \in p^{-1}(x_o)$. For such $x_o$, we shall denote the common
isotropy group of the points lying over $x_o$ by $\G_{x_o}$. Our second
assumption is that $X$ is connected, locally path connected, and locally
simply connected. These conditions are satisfied by all connected complex
algebraic varieties. Finally, we shall assume that if $g\in \G$ acts trivially
in the neighbourhood of some point of $X$, it acts trivially on all of $X$.
All three of these conditions are natural and will be satisfied in cases of
interest to us.

Let
$$
P(x_o) =
\{\text{homotopy classes of paths }([0,1],\{0,1\}) \to (X,p^{-1}(x_o))\}.
$$
Since $\G$ acts on the left of $X$, this is a left $\G$-set. Let
$$
Q(x_o) = \{(g,\gamma) \in \G \times P(x_o) : g^{-1}\cdot\gamma(0) = \gamma(1)\}.
$$
This has a natural left $\G$-action given by
$$
g : (h,\gamma) \mapsto (ghg^{-1},g\cdot \gamma).
$$

Define
$\piorb_1(\G\bs X,x_o)$ to be the quotient $\G\bs Q(x_o)$. This has a natural
group structure which can be understood by noting that $Q(x_o)$ is a groupoid.
The composition of two elements $(g,\gamma)$ and $(h,\mu)$ is defined when
$\gamma(1) = \mu(0)$. It is then given by
$$
(g,\gamma)\cdot (h,\mu) = (gh,\gamma\mu).
$$
To multiply two elements of $\piorb_1(G\bs X,x_o)$,
translate one of them until they are composable. For example, if we fix
a point $\xtilde_o$ of $p^{-1}(x_o)$, then each element of
$\piorb_1(G\bs X,x_o)$ has a representative of the form
$(g,\gamma)$ where $\gamma(0)=\xtilde_o$. This representation is unique up
to translation by an element of $\G_{x_o}$. To multiply two elements
$(g,\gamma)$ and $(h,\mu)$ starting at $\xtilde_o$, multiply $(h,\mu)$ by
$g$ so that it can be composed with $(g,\gamma)$. Then the product of these
two elements in $\piorb_1(G\bs X,x_o)$ is represented by the path
$$
(g,\gamma)\cdot(g^{-1}hg,g^{-1}\cdot \mu) = (hg, \gamma\, (g^{-1}\cdot \mu)).
$$

\begin{xca}
Show that if $k\in \G_{x_o}$, then $k$ acts 
trivially on $P(x_o)$. Deduce that the multiplication above is well defined.
(This is where we need $\G_{x_o}$ to be central in $\G$.)
\end{xca}

The following exercises should help give some understanding of orbifold
fundamental groups.

\begin{xca}
Show that if $\G$ acts fixed point freely on $X$, then there is a natural
isomorphism between $\piorb_1(\G\bs X,x_o)$ and the usual fundamental group
$\pi_1(\G\bs X,x_o)$. In particular, if we view a topological space $X$ as an
orbifold, then the two notions of fundamental group agree.
\end{xca}

\begin{xca}
Show that if $\G$ is an abelian group that acts trivially on a one point
space $X=\{*\}$, then $\piorb_1(\G\bs X,x_o)$ is defined and there is a
natural isomorphism
$$
\piorb_1(\G\bs X,x_o) \cong \G.
$$
\end{xca}

\begin{xca}
Let $X = \C$ and $\G = \Z/n\Z$. Define an action $\G$ on $X$ by
letting the generator 1 of $\Z/n\Z$ act by multiplication by $e^{2\pi i/n}$.
Show that $\piorb_1(\G\bs X,x_o)$ is defined for all $x_o \in \G\bs \C$, and
for each $x_o$ there is a natural isomorphism
$$
\piorb_1(\G\bs X,x_o) \cong \Z/n\Z.
$$
\end{xca}

\begin{xca}
Show that for each choice of a point $x\in p^{-1}(x_o)$, there is a
natural isomorphism
$$
\phi_x : \piorb_1(\G\bs X,x_o) \to \G
$$
defined by $\phi_x(g,\gamma) = g^{-1}$ when $\gamma(0) = x$.
Show that if $y = hx$, then $\phi_y = h\phi_x h^{-1}$.
\end{xca}

\begin{remark}
One can also define the orbifold fundamental group of $\G\bs X$ using the
Borel construction. First, find a contractible space $E\G$ on which $\G$
acts discontinuously and fixed point freely (any one will do). There is a
canonical construction of such spaces.
(See for example \cite[p.~19]{brown}.) Fix a base point $e_o \in E\G$. The
diagonal action of $\G$ on $E\G \times X$, a simply connected space, is fixed
point free and the map
$$
q:E\G \times X \to \G \bs (E\G\times X)
$$
is a covering mapping with Galois group $\G$. For $x \in X$, define
$$
\piorb_1(\G\bs X,q(e_o,x)) = \pi_1(\G(E\G\times X),q(e_o,x)).
$$
If $\G_x$ is central in $\G$, then this depends only on $x_o = p(x)$. It is
not difficult to show that, in this case, this definition agrees with the more
elementary one given above.
\end{remark}

The moduli space $\M_1$ is naturally an orbifold, being the quotient
of $\H$ by $SL_2(\Z)$. The condition on $x_o$ is satisfied for all
points not in the orbit of $i$ or $\rho$. Consequently, $\piorb_1(\M_1,x_o)$
is defined for all $x_o$ other than those corresponding to the orbits
of $i$ and $\rho$, and for each such $x_o$, there is an isomorphism
$$
\piorb_1(\M_1,x_o) \cong SL_2(\Z)
$$
which is well defined up to conjugacy.

This gives the following restatement of Proposition~\ref{class_map}.

\begin{proposition}
If $f : X \to T$ is a holomorphic family of genus 1 curves, then the
period mapping $\phi_f : T \to \M_1$ is a morphism of orbifolds.
\end{proposition}

\begin{remark}
\label{why_not}
Since $\M_1$ can also be written as $PSL_2(\Z)\bs \H$, it is natural
to ask why we are not giving $\M_1$ this orbifold structure. This question
will be fully answered in the section on the universal elliptic curve and
in subsequent sections on curves of higher genus. For the time being, just
note that if we give $\M_1$ the orbifold structure $PSL_2(\Z)\bs \H$, then
$\piorb(\M_1)$ would be $PSL_2(\Z)$ instead of $SL_2(\Z)$.
\end{remark}

\subsection{The Universal Elliptic Curve}

Let's attempt to construct a ``universal elliptic curve'' over $\M_1$.
We begin by constructing one over $\H$. The group $\Z^2$ acts on
$\C \times \H$ by
$$
(n,m) : (z,\tau) \mapsto (z+ n\tau + m, \tau).
$$
This action is fixed point free, so the quotient $\Z^2\bs(\C\times \H)$
is a complex manifold. The fiber of the projection
$$
\Z^2\bs(\C\times \H) \to \H
$$
over $\tau$ is simply the elliptic curve $\C/\Z \oplus \Z\tau$. The
family has the section that takes $\tau \in \H$ to the identity element
of the fiber lying over it. So this really is a family of elliptic curves.

Let's see what happens if we try to quotient out by $\Z^2$ and $SL_2(\Z)$
at the same time. First note that $SL_2(\Z)$ acts on $\Z^2$ on the right
by matrix multiplication.
We can thus form the semi-direct product $SL_2(\Z)\ltimes \Z^2$. This
is the group whose underlying set is $SL_2(\Z)\times \Z^2$ and whose
multiplication is given by
$$
(A_1,(n_1,m_1))(A_1,(n_2,m_2)) = (A_1A_2,(n_1,m_1)A_2 + (n_2,m_2)).
$$

\begin{xca}
Show that the action of $SL_2(\Z)\ltimes \Z^2$ on $\C \times \H$ given by
$$
\left(\begin{pmatrix} a & b \cr c & d \end{pmatrix},(n,m)\right) : (z,\tau)
\mapsto
\left(\frac{z + n\tau + m}{c\tau + d},\frac{a\tau + b}{c\tau + b}\right)
$$
is well defined. Show that there is a well defined projection
$$
(SL_2(\Z)\ltimes \Z^2) \bs (\C \times \H) \to SL_2(\Z) \bs \H.
$$
\end{xca}

This is a reasonable candidate for the universal curve. But we should
be careful.

\begin{xca}
Show that the fiber of the natural projection
\begin{equation}
\label{projn}
(SL_2(\Z)\ltimes \Z^2) \bs (\C \times \H) \to SL_2(\Z) \bs \H
\end{equation}
over the point corresponding to the elliptic curve $(C;0)$ is the quotient
of $C$ by the finite group $\Aut(C;0)$. In particular, the fiber of the 
projection over $[(C;0)]$ is always a quotient of $C/\pm \cong \P^1$, and
is never $C$.
\end{xca}

However, it is more natural to regard $\E:=(SL_2(\Z)\ltimes \Z^2)\ltimes\H$
as an orbifold. We shall do this. First note that the projection (\ref{projn})
has an orbifold section that is induced by the mapping
$$
\H \to \C \times \H
$$
that takes $\tau$ to $(0,\tau)$. Note that if we consider $\M_1$ as
the orbifold $PSL_2(\Z)\bs \H$, then the section would not exist. It is
for this reason that we give $\M_1$ the orbifold structure with orbifold
fundamental group $SL_2(\Z)$ instead of $PSL_2(\Z)$.

The following should illustrate why it is more natural to view $\M_1$ 
as an orbifold than as a variety.

\begin{theorem}
There is a natural one-to-one correspondence between holomorphic orbifold
mappings from a smooth complex curve (or variety) $T$ to $\M_1$ and families
of elliptic curves over $T$. The correspondence is given by pullback. \qed
\end{theorem}

We could define a `geometric points' of an orbifold $X$ to be an orbifold
map from a point with trivial fundamental group to $X$. For example, suppose
$\tau$ is any point of $\H$. Denote its isotropy group in $SL_2(\Z)$ by
$\G_\tau$. Then $\G_\tau\bs \{\tau\}$ is a point of $SL_2(\Z)\bs \H$, but is
not a geometric point as it has fundamental group $\G_\tau$, which is always
non-trivial as it contains $\Z/2\Z$. The corresponding geometric point
corresponds to the `universal covering' $\{\tau\} \to \G_\tau\bs \{\tau\}$
of this point. One can define pullbacks of orbifolds in such a way that the
fiber of the pullback of the universal curve to the geometric point $\tau$
of $\M_1$ is the corresponding elliptic curve $\C/(\Z \oplus \Z\tau)$.

Denote the $n$-fold fibered product of $\E \to \M_1$ with itself by
$\E^{(n)} \to \M_1$. It is naturally an orbifold (exercise). This has divisors
$\D_{jk}$ that consist of the points of $\E^{(n)}$ where the $j$th and $k$th
points agree. It also has divisors $\D_j$ where the $j$th point is zero.

\begin{xca}
Show that the moduli space $\M_{1,n}$ can be identified with
$$
\E^{(n-1)} - \big(\bigcup_{j<k} \D_{jk}\, \cup\ \bigcup_j \D_j\big).
$$
\end{xca}

\subsection{Modular Forms}

There is a natural orbifold line bundle $\cL$ over $\M_1$. It is the 
quotient of $\H \times \C$ by $SL_2(\Z)$ where $SL_2(\Z)$ acts via
the formula
$$
\begin{pmatrix} a & b \cr c & d \end{pmatrix} : (\tau,z) \to
\bigg(\frac{a\tau + b}{c\tau + d}, (c\tau + d)z\bigg).
$$

\begin{xca}
Show that this is indeed an action. Show that the $k$th power
$\cL^{\otimes k}$ of this line bundle is the quotient of $\H \times \C$
by the action
$$
\begin{pmatrix} a & b \cr c & d \end{pmatrix} : (\tau,z) \to
\bigg(\frac{a\tau + b}{c\tau + d}, (c\tau + d)^kz\bigg).
$$
\end{xca}

Orbifold sections of $\cL^{\otimes k}$ correspond to holomorphic
functions $f : \H \to \C$ for which the mapping
$$
\H \to \H \times \C \qquad \tau \mapsto (\tau,f(\tau))
$$
is $SL_2(\Z)$-equivariant.

\begin{xca}
Show that $f : \H \to \C$ corresponds to a section of $\cL^{\otimes k}$
if and only if
$$
f\bigg(
\begin{pmatrix} a & b \cr c & d \end{pmatrix}
\tau\bigg) = (c\tau + d)^k f(\tau).
$$
\end{xca}

Such a  function is called a {\it modular form of weight} $k$ for $SL_2(\Z)$.

\begin{xca}
Show that $f: \H \to \C$ is a modular form of weight $2k$ if and only if the
$k$-differential (i.e., section of the $k$th power of the canonical bundle)
$$
f(\tau) d\tau^{\otimes k}
$$
is invariant under $SL_2(\Z)$. Deduce that the (orbifold) canonical bundle of
$\M_1$ is isomorphic to $\cL^{\otimes 2}$.
\end{xca}

Some basic properties and applications of modular forms are given in
Chapter~VII of \cite{serre}.

Just as one can define the Picard group of a complex analytic variety to
be the group of isomorphism classes of holomorphic line bundles over it,
one can define the Picard group $\Picorb X$ of a holomorphic orbifold $X$ as
the group of isomorphism classes of orbifold line bundles over $X$. The
details are left as a straightforward exercise. The following result will
be proved later in the lectures.

\begin{theorem}
The Picard group of the orbifold $\M_1$ is cyclic of order $12$ and generated
by the class of $\cL$.
\end{theorem}

It is easy to see that the class of $\cL$ in $\Picorb \M_1$ as
$\cL^{\otimes 12}$ is trivialized by the cusp form
$$
\D(q) = (2\pi)^{12}q \prod_{n=1}^\infty (1-q^n)^{24}
$$
of weight 12 which has no zeros or poles in $\H$. Also, the non-existence
of meromorphic modular forms of weight $k$ with $0<k<12$ and no zeros or
poles in $\H$
implies that $\cL$ is an element of order 12 in $\Pic_\orb \M_1$. We will
show later that the order of $\Pic_\orb\M_1$ is of order 12, from which it
will follow that $\Pic_\orb \M_1$ is isomorphic to $\Z/12$ and generated by
the class of $\cL$.

\begin{center}
{\sc Lecture 2: Teichm\"uller Theory}
\end{center}

We have seen in genus 1 case that $\M_1$ is the quotient
$$
\G_1 \bs \X_1
$$
of a contractible complex manifold $\X_1 = \H$ by a discrete group
$\G_1 = SL_2(\Z)$. The action of $\G_1$ on $\X_1$ is said to be {\it virtually
free} --- that is, $\G_1$ has a finite index subgroup which acts (fixed point)
freely on $\X_1$.\footnote{More generally, we say that a group $G$ has a
property $P$ {\it virtually} if $P$ holds for some finite index subgroup of
$G$.} In this section we will generalize this to all $g\ge 1$ --- we will
sketch a proof
that there is a contractible complex manifold $\X_g$, called {\it Teichm\"uller
space}, and a group $\G_g$, called the {\it mapping class group}, which acts
virtually freely on $\X_g$. The moduli space of genus $g$ compact Riemann
surfaces is the quotient:
$$
\M_g = \G_g \bs \X_g.
$$
This will imply that $\M_g$ has the structure of a complex analytic variety
with finite quotient singularities.

Teichm\"uller theory is a difficult and technical subject. Because of this,
it is only possible to give an overview.

\section{The Uniformization Theorem}

Our basic tool is the the generalization of the Riemann Mapping Theorem
known as the Uniformization Theorem. You can find a proof, for example, in
\cite{forster} and \cite{farkas-kra}.

\begin{theorem}
Every simply connected Riemann surface\footnote{We follow the convention
that every Riemann surface is, by definition, connected} is biholomorphic
to either $\P^1$, $\C$ or $\H$.
\end{theorem}

\begin{xca}
Show that no two of $\P^1$, $\H$ and $\C$ are isomorphic as Riemann surfaces.
\end{xca}

We have already seen that $\Aut \P^1 \cong PSL_2(\C)$ and that
$PSL_2(\R) \subseteq \Aut \H$.

\begin{xca}
\begin{enumerate}
\item Show that
$$
\Aut \C = \{z \mapsto az + b : a \in \C^\ast \text{ and } b\ \in \C \}.
$$
\item Show that there is an element $T$ of $\Aut \P^1$ that restricts to
an isomorphism $T : \D \to \H$ between the unit disk and the upper half
plane.
\item Show that every element of $\Aut \D$ is a fractional linear
transformation. Hint: Use the Schwartz Lemma.
\item Deduce that every element of $\Aut \H$ is a fractional linear
transformation and that $\Aut \H \cong PSL_2(\R)$.
\end{enumerate}
\end{xca}

If $X$ is a Riemann surface and $x\in X$, then $\pi_1(X,x)$ acts on
the universal covering $\Xtilde$ as a group of biholomorphisms. This
action is fixed point free.

\begin{xca}
\label{quotient}
Show that if $X$ is a Riemann surface whose universal covering is isomorphic
to $\H$ and $x\in X$, then there is a natural injective homomorphism
$\rho : \pi_1(X,x) \to PSL_2(\R)$ which is injective and has discrete image.
Show that this homomorphism is unique up to conjugation by an element of
$PSL_2(\R)$, and that the conjugacy class of $\rho$ is independent
of the choice of $x$. Show that $X$ is isomorphic to $\im \rho \bs \H$ and
that the conjugacy class of $\rho$ determines $X$ up to isomorphism.
\end{xca}

This will give a direct method of putting a topology on $\M_g$ when
$g\ge 2$. But first some preliminaries.

\begin{xca}
Show that if $X$ is a Riemann surface whose universal covering is
\begin{enumerate}
\item $\P^1$, then $X$ is isomorphic to $\P^1$;
\item $\C$, then $X$ is isomorphic to $\C$, $\C^\ast$ or is a genus 1 curve.
\end{enumerate}
Hint: Classify the subgroups of $\Aut \P^1$ and $\Aut \C$ that act properly
discontinuously and freely.
\end{xca}

So we come to the striking conclusion that the universal covering of a
Riemann surface not isomorphic to $\P^1$, $\C$, $\C^\ast$ or a genus 1
curves must be isomorphic to $\H$. In particular, the universal covering
of $\C - \{0,1\}$ is $\H$. Picard's Little Theorem is a consequence.
 
\begin{xca}
Show that if $f: \C \to \C$ is a holomorphic mapping
that omits 2 distinct points of $\C$, then $f$ is constant.
\end{xca}

For our immediate purposes, the most important fact is that the universal
covering of every
Riemann surface of genus $g\ge 2$ is isomorphic to $\H$. More generally,
we have the following interpretation of the stability condition
$2g-2+n>0$.

\begin{xca}
Show that if $U$ is a Riemann surface of the form $X - F$ where $X$ is a
compact Riemann surface genus $g$ and $F$ a finite subset of $X$ of cardinality
$n$, then the universal covering of $U$ is isomorphic to $\H$ if and only if
$2g-2 + n > 0$. Observe that $2-2g-n$ is the topological Euler characteristic
of $U$.
\end{xca}

\section{Teichm\"uller Space}

Suppose that $g\ge 2$ and that $S$ is a compact oriented surface of genus
$g\ge 2$. Fix a base point $x_o \in S$ and set $\pi = \pi_1(S,x_o)$. Let
$$
\X_g = \left\{\parbox{3.25in}{conjugacy classes of injective representations
$\rho : \pi \hookrightarrow PSL_2(\R)$ such that $\im \rho$ acts freely on
$\H$ and $\im\rho\bs \H$ is a compact Riemann surface of genus $g$} \right\}
$$

It is standard that we can choose generators $a_1,\dots,a_g$, $b_1,\dots,b_g$
of $\pi$ that are subject to the relation
$$
\prod_{j=1}^g [a_j,b_j] = 1.
$$
where $[x,y]$ denotes the commutator $xyx^{-1}y^{-1}$.

\begin{xca}
Show that $\X_g$ can also be identified naturally with the set
$$
\left\{\parbox{4.25in}{compact Riemann surfaces $C$ plus a conjugacy class
of isomorphisms $\pi \cong \pi_1(C)$ modulo the obvious isomorphisms} \right\}
$$
\end{xca}

The group $PSL_2(\R)$ is a real affine algebraic group. It can be realized as
a closed subgroup of $GL_3(\R)$ defined by real polynomial equations.

\begin{xca}
Denote the standard 2-dimensional representation of $SL_2(\R)$ by $V$. Denote
the second symmetric power of $V$ by $S^2V$. Show that this is a 3-dimensional
representation of $SL_2(\R)$. Show that $-I\in SL_2(\R)$ acts trivially on
$S^2V$, so that $S^2 V$ is a 3-dimensional representation of $PSL_2(\R)$.
Show that this representation is faithful and that its image is defined by
polynomial equations. Deduce that $PSL_2(\R)$ is a real affine algebraic
group.
\end{xca}

Because of this, we will think of elements of $PSL_2(\R)$ as $3\times 3$
matrices whose entries satisfy certain polynomial equations.
A representation $\rho : \pi \to PSL_2(\R)$ thus corresponds to a
collection of matrices
$$
A_1,\dots, A_g,B_1,\dots, B_g
$$
in $PSL_2(\R) \subset GL_3(\R)$ that satisfy the polynomial equation
\begin{equation}
\label{reln}
I - \prod_{j=1}^g [A_j,B_j] = 0;
\end{equation}
the correspondence is given by setting $A_j = \rho(a_j)$ and $B_j = \rho(b_j)$.

The set of all representations $\rho : \pi \to PSL_2(\R)$ is the real
algebraic subvariety $\cR$ of $PSL_2(\R)^{2g}$ consisting of all $2g$-tuples
$$
(A_1,\dots,A_g,B_1,\dots,B_g)
$$
that satisfy (\ref{reln}). This is a closed subvariety of $PSL_2(\R)^{2g}$
and is therefore an affine variety. Note that\footnote{One has to be  more
careful here as, in real algebraic geometry, affine varieties of higher
codimension can all be cut out by just one equation. Probably the best way
to proceed is to compute the derivative of the product of commutators
map $PSL_2(\R)^{2g} \to PSL_2(\R)$ at points in the fiber over the 
identity.}
$$
\dim \cR \ge 2g\cdot \dim PSL_2(\R) - \dim PSL_2(\R) = 6g-3.
$$
Observe that $PSL_2(\R)$ acts on $\cR$
on the right by conjugation:
$$
A : \rho \mapsto \{g \mapsto A^{-1}\rho(g)A\}.
$$

\begin{xca}
Show that
$$
U := \left\{
\parbox{2.75in}{$\rho :\pi \hookrightarrow PSL_2(\R)$ such that $\im\rho$ acts
freely on $\H$ and $\im\rho\bs\H$ is a compact Riemann surface of genus $g$}
\right\}
$$
is an open subset of $\cR$ and that it is closed under the $PSL_2(\R)$-action.
Hint: One can construct a fundamental domain $F$ for the action of $\pi$ on
$\H$ given by $\rho$ as follows. First choose a point $x\in \H$. Then
take $F$ to be all points in $\H$ that are closer to $x$ than to any of
the points $\rho(g)x$ (with respect to hyperbolic distance) where $g\in \pi$
is non-trivial. Then $F$ is a compact subset of $\H$ whose boundary is a
union of geodesic segments. As $\rho$ varies continuously, the orbit of
$x$ varies continuously. Now choose $x$ to be generic enough and study the
change in $F$ as $\rho$ varies.
\end{xca}

Note that $\X_g = U/PSL_2(\R)$.

\begin{xca}
Show that the center of $PSL_2(\R)$ is trivial.
\end{xca}

\begin{proposition}
\label{no_inv}
If $\rho \in U$, then the stabilizer in $PSL_2(\R)$ of $\rho$ is trivial.
\end{proposition}

\begin{proof}
The first step is to show that if $\rho \in U$, then $\im\rho$ is Zariski
dense in $PSL_2(\R)$. To prove this, it suffices to show that $\im\rho$ is
Zariski dense in the set of complex points $PSL_2(\C)$. From Lie theory
(or algebraic group theory), we know that all proper subgroups of $PSL_2(\C)$
are extensions of a finite group by a solvable group. Since $\im \rho$ is
isomorphic to $\pi$, and since $\pi$ is not solvable (as it contains a free
group of rank $g$), $\im\rho$ cannot be contained in any proper algebraic
subgroup of $PSL_2(\R)$, and is therefore Zariski dense.

An element $A$ of $PSL_2(\R)$ stabilizes $\rho$ if and only if
$A^{-1}\rho(g) A = \rho(g)$ for all $g\in\pi$ --- that is, if and only if
$A$ centralizes $\im\rho$. But $A$ centralizes $\im\rho$ if and only if
it centralizes the Zariski closure of $\im\rho$. Since $\im\rho$ is Zariski
dense in $PSL_2(\R)$, and since the center of $PSL_2(\R)$ is trivial, it
follows that $A$ is trivial.
\end{proof}

\begin{theorem}
The set $U$ is a smooth manifold of real dimension $6g-3$. The group
$PGL_2(\R)$ acts principally on $U$, so that quotient $\X_g = U/PSL_2(\R)$
is a manifold of real dimension $6g-6$. 
\end{theorem}

\begin{proof}[Sketch of Proof.]
It can be shown, using deformation theory, that the Zariski tangent
space of $\cR$ at the representation $\rho$ of $\pi$ is given by the
relative cohomology group:\footnote{This is a non-standard way of giving
this group --- Weil used cocycles. Here $A^\dot(S;\bA)$ is the complex of
differential forms on $S$ with coefficients in the local system $\bA$. It
is a differential graded Lie algebra. Restricting to the base point $x_o$
one obtains an augmentation $A^\dot(S;\bA) \to \sl_2(\R)$.
The complex $A^\dot(S,x_o;\bA)$ is the kernel of this, and $H^\dot(S,x_o;\bA)$
is defined to be the cohomology of this complex. Briefly, this is related to
deformations of $\rho$ as follows: when the representation $\rho$ is deformed,
the vector bundle underlying the local system corresponding to $\rho$ does
not change, only the connection on it. But a new connection on the bundle
corresponding to $\rho$ will differ from the original connection by an
element of $A^1(S;\bA)$. So a deformation of $\rho$ corresponds to a family
of sections $w(t)$ of $A^1(S;\bA)$ where $w(0) = 0$. For each $t$ the new 
connection is flat. This will imply that the $t$ derivative of $w(t)$ at
$t=0$ is closed in $A^\dot(S,\bA)$. Changes of gauge that do not alter the
marking of the fiber over $x_o$ are, to first order, elements of
$A^0(S,x_o;\bA)$. The $t$ derivative of such first order changes of gauge
gives the trivial first order deformations of $\rho$.}
$$
T_\rho \cR = H^1(S,x_o;\bA)
$$
where $\bA$ is the local system (i.e., locally constant sheaf) over $S$ whose
fiber over $x_o$ is $\sl_1(\R)$ and whose monodromy representation is
the homomorphism
$$
\begin{CD}
\pi @>\rho>> PSL_2(\R) @>{\Ad}>> \Aut \sl_2(\R).
\end{CD}
$$
Here $\sl_2(\R)$ denotes the Lie algebra of $PSL_2(\R)$, which is the
set of $2\times 2$ matrices of trace 0, and $\Ad$ is the adjoint representation
$$
A \mapsto \{X \mapsto AXA^{-1}\}.
$$
A result equivalent to this was first proved by Andr\'e Weil \cite{weil}.
The long exact sequence of the pair $(S,x_o)$ with coefficients in $\bA$
gives a short exact sequence
$$
0 \to \sl_2(\R) \to H^1(S,x_o;\bA) \to H^1(S;\bA) \to 0
$$
and an isomorphism $H^2(S,x_o;\bA) \cong H^2(S;\bA)$.

Deformation theory also tells us that if $H^2(S,x_o;\bA)$ vanishes,
then $\cR$ is smooth at $[\rho] \in \cR$.\footnote{A proof of this and also a
proof of Weil's result can be found in the beautiful paper of Goldman and
Millson \cite{goldman-millson}.} The Killing form is the symmetric bilinear
form on $\sl_2(\R)$ given by
$$
b(X,Y) = \tr(XY).
$$
It is non-degenerate and invariant under $\Ad\circ\rho$. It follows that $\bA$
is isomorphic to its dual $\bA^\ast$ as a local system. Thus it follows from
Poinar\'e duality (with twisted coefficients) that
$$
H^2(S,\bA) \cong H^0(S,\bA)^\ast.
$$
But if $\rho$ is an element of $U$, it follows from Proposition~\ref{no_inv}
that $H^0(S,\bA)$ vanishes. Consequently, if $\rho \in U$,
then
$$
\dim H^1(S,\bA) = - \chi(S,\bA) = -(\rank\bA)\cdot \chi(S) = 3(2g-2) = 6g - 6
$$
and
$$
\dim H^1(S,x_o;\bA) = 6g-3.
$$
We saw above that $U$ has dimension $\ge 6g-3$, but $\dim U$ is not bigger
than the dimension of any of its Zariski tangent spaces, which we have
just shown is $6g-3$. It follows that $\dim U = 6g-3$ and that it is smooth
as its dimension equals the rank of its Zariski tangent space. Finally,
since the $PSL_2(\R)$-action on $U$ is principal, $\X_g = U/PSL_2(\R)$ is
smooth of dimension $(6g-3) - \dim PSL_2(\R) = 6g-6$.
\end{proof}

It is not very difficult to show that the tangent space of $\X_g$ at $[\rho]$
is naturally isomorphic to $H^1(S,\bA)$.

\section{Mapping Class Groups}

The mapping class group $\G_g$ plays a role for higher genus ($g\ge 2$)
Riemann surfaces analogous to that played by $SL_2(\Z)$ in the theory of
genus 1 curves. As in the previous section, $S$ is a fixed compact oriented
surface of genus $g$, $x_o \in S$, and $\pi = \pi_1(S,x_o)$. But unlike
the last section, we allow any $g\ge 1$.

The mapping class group $\G_S$ of $S$ is defined to be the group of connected
components of the group of orientation preserving diffeomorphisms of $S$:
$$
\G_S = \pi_0 \Diff^+ S.
$$
Here $\Diff S$ is given the compact open topology. It is the finest
topology on $\Diff S$ such that a function $f : K \to \Diff S$
from a compact space into $\Diff S$ is continuous if and only if the 
map $\phi_f : K \times S \to S$ defined by
$$
\phi_f (k,x) = f(k)(x)
$$
is continuous.

\begin{xca}
Suppose that $M$ is a smooth manifold.
Show that a path in $\Diff M$ from $\phi$ to $\psi$ is a homotopy
$\Psi : [0,1]\times M \to M$ from $\phi$ to $\psi$ such that for each
$t \in [0,1]$, the function $\Psi(t,\blank) : M \to M$ that takes $x$ to
$\Psi(t,x)$ is a diffeomorphism. Such a homotopy is called an {\it isotopy}
between $\phi$ and $\psi$.
\end{xca}

Recall that the outer automorphism group of a group $G$ is defined to
be the quotient
$$
\Out G = \Aut G/ \Inn G
$$
of the automorphism group $\Aut G$ of $G$ by the subgroup $\Inn G$ of inner 
automorphisms.

\begin{xca}
Suppose that $M$ is a smooth manifold. Show that each diffeomorphism
of $M$ induces an outer automorphism of its fundamental group and
that the corresponding function
$$
\Diff M \to \Out \pi(M)
$$
is a homomorphism. (Hint: Fix a base point of $M$ and show that every
diffeomorphism of $M$ is isotopic to one that fixes the base point.)
Show that the kernel of this homomorphism contains
the subgroup of all diffeomorphisms isotopic to the identity. Deduce that
there is a homomorphism
$$
\pi_0(\Diff M) \to \Out \pi_1(M).
$$
\end{xca}

Taking $M$ to be $S$, our reference surface, we see that there is a 
natural homomorphism
$$
\Phi : \G_S \to \Out \pi.
$$
It is a remarkable fact, due to Dehn \cite{dehn} and Nielsen \cite{nielsen},
that this map is almost an isomorphism. (I believe there is a proof of this
in the translation of some papers of Dehn on topology by Stillwell.)

\begin{theorem}[Dehn-Nielsen]
For all $g \ge 1$,
the homomorphism $\Phi$ is injective and the image of $\Phi$ is a subgroup
of index $2$. \qed
\end{theorem}

One can describe the index 2 subgroup of $\Out \pi$ as follows. Since
$S$ has genus $\ge 1$, the universal covering of $S$ is contractible, and
there is a natural isomorphism $H_i(S,\Z) \cong H_i(\pi,\Z)$. Each element
of $\Out \pi$ induces an automorphism of $H_i(\pi,\Z)$ as inner automorphisms
of a group act trivially on its cohomology. The image of $\Phi$ is the kernel
of the natural homomorphism
$$
\Out \pi \to \Aut H_2(\pi,\Z) \cong \{\pm 1\}.
$$

When $g=1$, $\pi \cong \Z^2$, and
$$
\Out \pi = \Aut \pi \cong GL_2(\Z).
$$

\begin{xca}
Show that when $g=1$, the homomorphism $\Out \pi \to \{\pm 1\}$ corresponds
to the determinant $GL_2(\Z) \to \{\pm 1\}$.
\end{xca}

The genus 1 case of the Dehn-Nielsen Theorem can be proved by elementary
means:

\begin{xca}
Show that if $g=1$, then the homomorphism $\G_S \to SL_2(\Z)$ is an
isomorphism. Hint: Construct a homomorphism from $SL_2(\Z)$ to $\Diff S$.
\end{xca}

Since two compact orientable surfaces are diffeomorphic if and only if
they have the same genus, it follows that the group $\G_S$ depends only
on the genus of $S$. For this reason, we define $\G_g$ to be $\G_S$ where
$S$ is any compact genus $g$ surface.

\section{The Moduli Space}

In this section, $S$ is once again a compact oriented surface of genus
$g\ge 2$ and $\pi$ denotes its fundamental group. The mapping class
group $\G_g = \G_S$ acts smoothly on Teichm\"uller space $\X_g$ on the left
via the homomorphism $\G_g \to \Out \pi$.

\begin{xca}
Describe this action explicitly. Show that the isotropy group of
$[\rho] \in \X_g$ is naturally isomorphic to the automorphism group of
the Riemann surface $\im\rho \bs \H$. (Compare this with the genus 1
case.) Deduce that all isotropy groups are finite.
\end{xca}

\begin{xca}
Show that two points in $\X_g$ are in the same $\G_g$ orbit if and only
if they determine the same Riemann surface. Hint: Use the Uniformization
Theorem.
\end{xca}

Rephrasing this, we get:

\begin{theorem}
If $g\ge 2$, then $\M_g$ is naturally isomorphic to the quotient of
$\X_g$ by $\G_g$. \qed
\end{theorem}

This result allows us to put a topology on $\M_g$. We give it the unique
topology so that $\X_g \to \M_g$ is a quotient mapping.
To try to understand the topology of $\X_g$ we shall need the following
fundamental result.

\begin{theorem}[Teichm\"uller]
For all $g\ge 2$, the Teichm\"uller space $\X_g$ is contractible and the
action of $\G_g$ on it is properly discontinuous.
\end{theorem}

A natural way to approach the proof of this theorem via hyperbolic geometry.
We do this in the next two sections.

\section{Hyperbolic Geometry}

Perhaps the most direct way to approach the study of Teichm\"uller space
is via hyperbolic geometry. The link between complex analysis and hyperbolic
geometry comes from the fact that the upper half plane has a complete metric
$$
ds^2 = \frac{1}{y^2}(dx^2 + dy^2)
$$
of constant curvature $-1$ whose group of orientation preserving isometries
is $PSL_2(\R)$. The Riemannian manifold $(\H,ds^2)$ is called the
{\it Poincar\'e upper half plane}. A good reference for hyperbolic geometry is
the book by Beardon \cite{beardon}.

\begin{xca}
Show that
\begin{enumerate}
\item every element of $PSL_2(\R)$ is an orientation preserving isometry of
$(\H,ds^2)$;
\item $PSL_2(\R)$ acts transitively on $\H$ and that the stabilizer of one
point (say $i$), and therefore all points, of $\H$ is isomorphic to $SO(2)$
and acts transitively on the unit tangent space of $\H$ at this point.
\end{enumerate}
Deduce that the Poincar\'e metric $ds^2$ has constant curvature and that
$PSL_2(\R)$ is the group of all orientation preserving isometries of the
Poincar\'e upper half plane.
\end{xca}

The fact that the group of biholomorphisms of $\H$ coincides with the group
of orientation preserving isometries of $\H$ is fundamental and is used as
follows. Suppose that $X$ is a Riemann surface whose universal covering is
$\H$. Then, by Exercise~\ref{quotient}, $X$ is biholomorphic to the quotient
of $\H$ by a discrete subgroup $\G$ of $PSL_2(\R)$ isomorphic to $\pi_1(X)$.
Now, since $PSL_2(\R)$ is also the group of orientation preserving isometries
of $\H$, the Poincar\'e metric $ds^2$ is invariant under $\G$ and therefore
descends to $X = \G\bs \H$. This metric is complete as the Poincar\'e metric
is.

\begin{xca}
Prove the converse of this: if $X$ is an oriented surface with a complete
hyperbolic metric, then $X$ has a natural complex structure whose canonical
orientation agrees with the original orientation.
\end{xca}

These two constructions are mutually inverse. Thus we conclude that when
$g \ge 2$, there is a natural one-to-one correspondence
\begin{multline}
\label{cplex-hyp}
\left\{
\parbox{2.2in}{isomorphism classes of compact Riemann surfaces
of genus $g$} \right\} \leftrightarrow \cr
\left\{
\parbox{3in}{isometry classes of compact oriented surfaces
of genus $g$ with a hyperbolic metric} \right\}
\end{multline}

It follows that if $g \ge 2$, then
$$
\M_g = \left\{
\parbox{3in}{isometry classes of compact oriented surfaces
of genus $g$ with a hyperbolic metric} \right\}.
$$

\begin{remark}
Likewise, a Riemann surface of genus 1 has a flat metric (unique up to
rescaling) which determines the complex structure. So $\M_1$ can be regarded
as the moduli space of (conformal classes of) flat tori.
\end{remark}

\begin{xca}
Show that if $X$ is a Riemann surface whose universal covering is $\H$,
then the group of orientation preserving isometries of $X$ equals the
group of biholomorphisms of $X$.
\end{xca}

\section{Fenchel-Nielsen Coordinates}

Once again, we assume that $g\ge 2$ and that $S$ is a reference surface of
genus $g$. The interpretation of Teichm\"uller space in terms of hyperbolic
geometry allows us to define coordinates on $\X_g$. To do this we decompose
the surface into ``pants.''

A {\it pair of pants} is a compact oriented surface of genus 0 with 3
boundary components. Alternatively, a pair of pants is a disk with 2 holes.
\begin{figure}[!ht]
\epsfig{file=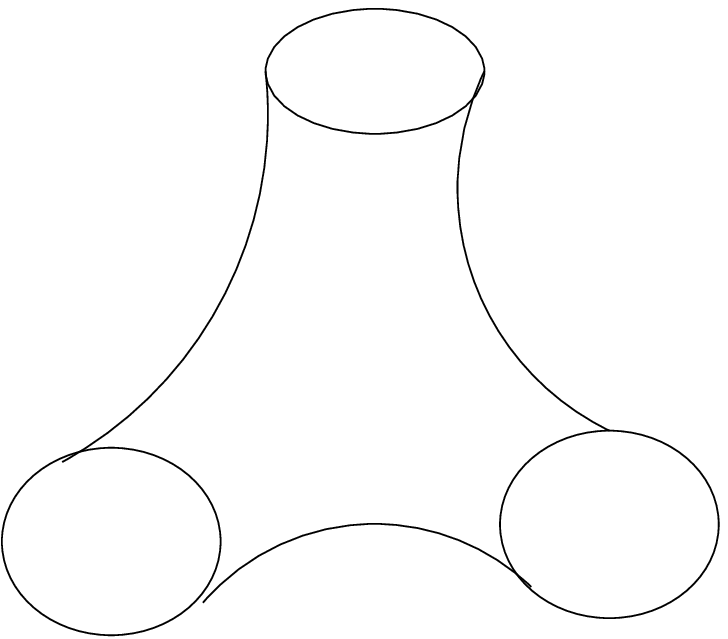, height=1in}
\hspace{1cm}
\epsfig{file=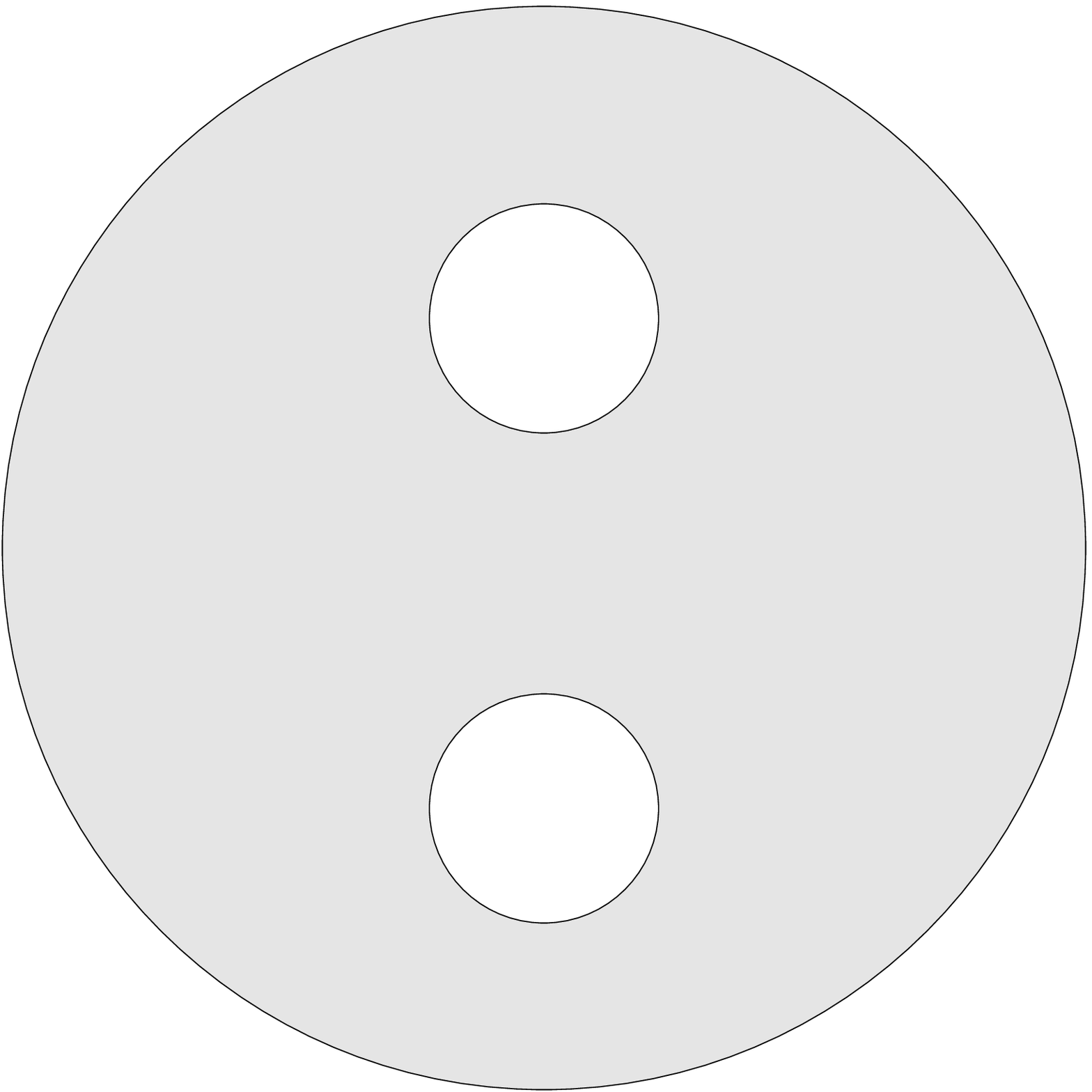, height=1in}
\caption{pairs of pants}
\end{figure}
A {\it simple closed curve} (SCC) on $S$ is an imbedded circle. A {\it
pants decomposition} of $S$ is a set of disjoint simple closed curves in $S$ 
that divides $S$ into pairs of pants.
\begin{figure}[!ht]
\epsfig{file=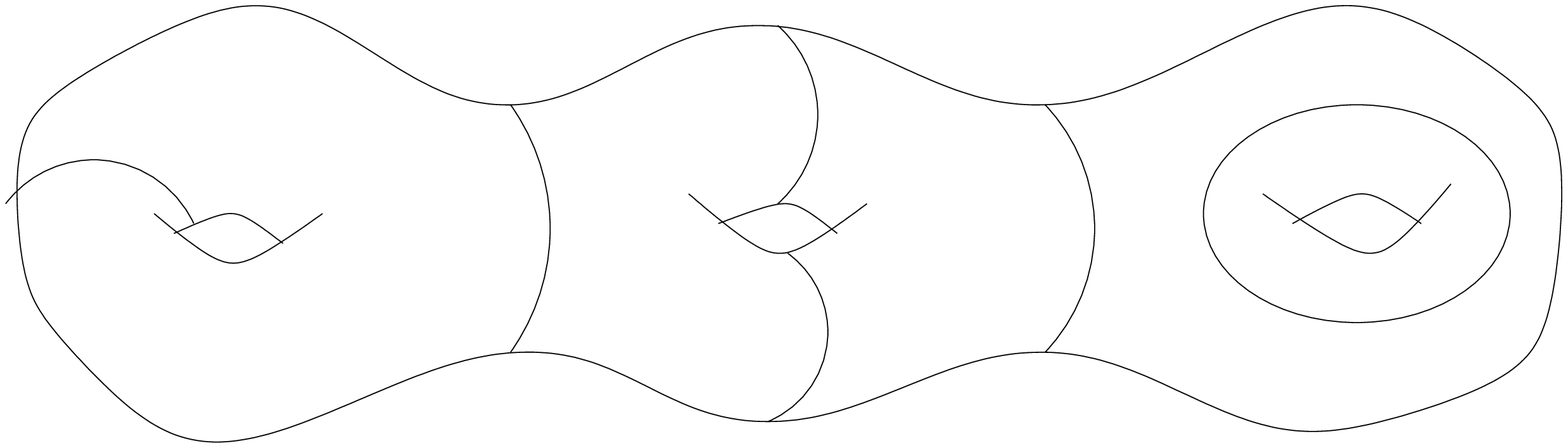, height=1in}
\caption{a pants decomposition}
\end{figure}

\begin{proposition}
If the simple closed curves $C_1,\dots,C_n$ divide $S$ into $m$ pairs of pants,
then $n=3g-3$ and $m=2g-2$.
\end{proposition}

\begin{proof}
Since a pair of pants has the homotopy type of a bouquet of 2 circles, it
has Euler characteristic $-1$. Since the Euler characteristic of a disjoint
union of circles is 0, we have
$$
2-2g = \chi(S) = m\cdot\chi(\text{a pair of pants}) = -m.
$$
Thus $m=2g-2$. Since each pair of pants has 3 boundary components, and
since each circle lies on the boundary of two pairs of pants, we see that
$$
n = 3m/2 = 3g-3.
$$
\end{proof}

If $T$ is a compact oriented surface with boundary and $\chi(T) < 0$, then
$T$ has a hyperbolic metric such that each boundary component is totally
geodesic. As in the case where the boundary of $T$ is empty, the set of
all hyperbolic structures on $T$ (modulo diffeomorphisms isotopic to the
identity) can be described as a representation variety. I will omit the 
details. This space is called the Teichm\"uller space of $T$ and will be
denoted by $\X_T$.

\begin{proposition}
If $P$ is a pair of pants, then $\X_P$ is a manifold diffeomorphic to $\R^3$.
The diffeomorphism is given by taking a hyperbolic structure on $P$ to the
lengths of its 3 boundary components.
\end{proposition}

\begin{proof}[Discussion of the Proof.]
A proof can be found in the Expos\'e 3, partie II by Po\'enaru in the book
\cite{flp}. The basic idea is that a pair of pants can be cut into two
isomorphic hyperbolic right hexagons, and that two right hyperbolic hexagons
are equivalent if the lengths of every other side of one equal the lengths
of the corresponding sides of the other. Also, one can construct hyperbolic
right hexagons where these lengths are arbitrary positive real numbers.
\end{proof}

A basic fact in hyperbolic
geometry is that if one has a compact hyperbolic surface $T$ (not assumed to
be connected) with totally geodesic boundary, then one can identify pairs
of boundary components of the same lengths to obtain a hyperbolic surface
whose new hyperbolic structure agrees with that on $T - \partial T$. Another
elementary fact we shall need is that if $C$ is a SCC in a compact hyperbolic
surface $T$ with totally geodesic boundary, then there is a unique closed
geodesic $\gamma : S^1 \to T$ that is freely homotopic to $C$.

Now suppose that $\cP=\{P_1,\dots,P_{2g-2}\}$ is a pants decomposition of
$S$ determined by the SCCs $C_1,\dots,C_{3g-3}$. For each hyperbolic structure
on $S$, we can find closed geodesics $\gamma_1,\dots, \gamma_{3g-3}$ that
are isotopic to the $C_j$. Taking the lengths of these, we obtain a function
$$
\ell : \X_S \to \R_+^{3g-3}.
$$
It is not difficult to show that this mapping is continuous. (See \cite{flp},
for example.) On the other hand, we can define an action of $\R^{3g-3}$ on
$\X_S$ and whose orbits lie in the fibers of $\ell$. To define the action
of $(\theta_1,\dots, \theta_{3g-3})$ on a hyperbolic surface $S$, write
$S$ as the union of $2g-2$ pairs of hyperbolic pants whose boundaries are
represented by geodesics $\gamma_1,\dots,\gamma_{3g-3}$ isotopic to the
$C_j$. Then twist the gluing map at $\gamma_j$ by an angle $\theta_j$.
This is also continuous.

\begin{theorem}[Douady,{\cite[Expos\'e 7]{flp}}]
If $g\ge 2$, then the map $\ell : \X_S \to \R_+^{3g-3}$ is a principal
$\R^{3g-3}$ bundle with the action described above. \qed
\end{theorem} 

\begin{corollary}
For all $g\ge 2$, Teichm\"uller space $\X_g$ is diffeomorphic to $\R^{6g-6}$.
\qed
\end{corollary}

\begin{xca}
At first it may appear that $\ell : \X_g \to \R_+^{3g-3}$ should be a
principal $(S^1)^{3g-3}$ bundle. Show that rotation by $2\pi$ about one
of the SCCs $C_j$ alters the representation $\rho$ by an automorphism
of $\pi$ that is not inner. Hint: it may help to first read the part
of Section~\ref{relations} on Dehn twists.
\end{xca}

\section{The Complex Structure}

Suppose that $C$ is a compact Riemann surface. Recall that a deformation of
$C$ is the germ about $t_o \in T$ of a proper analytic mapping $f:\cC \to T$,
where $T$ is an analytic variety, and an isomorphism $j : C \to f^{-1}(t_o)$.
One defines maps between deformations to be cartesian squares. A deformation
of $C$ is called universal if every other deformation is pulled back from it.
A standard result in deformation theory is that every Riemann surface of genus
$g\ge 2$ has a universal deformation. This has the property that $T$ is
smooth at $t_o$ and that $T$ is of complex dimension $3g-3$ with tangent
space at $t_o$ canonically isomorphic to $H^1(C,\Theta_C)$, where
$\Theta_C$ denotes the sheaf of holomorphic sections of the tangent bundle of
$C$.

As explained in Looijenga's lectures, $\M_g$ can be obtained by patching
such local deformation spaces together. If $C$ is a smooth projective curve
of genus $g$ and $\rho : \pi \to PSL_2(\R)$ is a representation such that
$C \cong \im \rho \bs \H$, then $[\rho] \in \X_g$ goes to $[C]\in \M_g$
under the projection $\X_g \to \M_g$. The following fact can be proved
using partial differential equations.


\begin{theorem}
If $(T,t_o)$ is a universal deformation space for $C$, then there is a
smooth mapping $\psi : T \to \X_g$ such that $\psi(t_o) = [\rho]$ which is
a diffeomorphism in a neighbourhood of $t_o$. Moreover, the composition of
$\psi$ with the projection $\X_g \to \M_g$ is the canonical mapping that
classifies the universal deformation of $C$. \qed
\end{theorem}

\begin{corollary}
There is a canonical $\G_g$ invariant complex structure on $\X_g$ such
that the projection $\X_g \to \M_g$ is a complex analytic mapping. \qed
\end{corollary}

It should be noted that, although $\X_g$ is a complex manifold diffeomorphic
to $\R^{6g-6}$, it is not biholomorphic to either a complex $3g-3$ ball or
$\C^{3g-3}$ when $g>1$.

Combining these results, we obtain the following basic result:

\begin{theorem}
The moduli space $\M_g$ is a complex analytic variety whose singularities
are all finite quotient singularities. Furthermore, $\M_g$ can be regarded
as an orbifold whose universal covering is the complex manifold $\X_g$ and
whose orbifold fundamental groups can be identified, up to conjugacy, with
$\G_g$. More precisely, if $C$ is a compact genus $g$ curve with no non-trivial
automorphisms, then there is a canonical isomorphism
$$
\piorb_1(\M_g,[C]) \cong \pi_0 \Diff^+ C. \qed
$$
\end{theorem}

\section{The Teichm\"uller Space $\X_{g,n}$}

It is natural to expect that there is a complex manifold $\X_{g,n}$ and
a discrete group $\G_{g,n}$ that acts holomorphically on $\X_{g,n}$ in such
a way that $\M_{g,n}$ is isomorphic to the quotient $\G_{g,n}\bs \X_{g,n}$.
We give a brief sketch of how this may be deduced from the results when
$n=0$ and $g\ge 2$.

First, we fix an $n$-pointed reference surface of genus $g$ where $g \ge 2$.
That is, we fix a compact oriented surface  $S$ and a subset
$P=\{x_1,\dots,x_n\}$ of $n$ distinct points of $S$. Define
$$
\G_{g,n} := \pi_0 \Diff^+(S,P)
$$ 
By definition, elements of $\Diff^+(S,P)$ are orientation preserving and
act trivially on $P$.

Here is a sketch of a construction of $\X_{g,1}$. One can construct the
$\X_{g,n}$ when $n>1$ in a similar fashion.

There is a universal curve $\cC \to \X_g$. This can be constructed
using deformation theory. This is a fiber bundle in the topological sense
and it is not difficult to see that the action of the mapping class group
$\G_g$ on $\X_g$ can be lifted to an action on $\cC$ such that the projection
is equivariant. Since the fiber of $\cC \to \X_g$ is a compact surface of
genus $g$, and since $\X_g$ is contractible, $\cC$ has the homotopy type
of a surface of genus $g$ and therefore has fundamental group isomorphic to
$\pi_1(S)$. Define $\X_{g,1}$ to be the universal covering of $\cC$. This is
a complex manifold as $\cC$ is. The fibers of the projection $\X_{g,1}$ are
all isomorphic to $\H$. Since $\X_g$ is contractible, so is $\X_{g,1}$.

\begin{xca}
Show that there is a natural bijection
$$
\X_{g,1} = \left\{\parbox{3.25in}{conjugacy classes of representations
$\rho:\pi \to PSL_2(\R)$ such that $C_\rho := \im\rho \bs \H$ is of genus
$g$, plus a point $x\in C_\rho$} \right\}.
$$
Use this (or otherwise) to show that $\G_{g,1}$ acts on $\X_{g,1}$ and
that the quotient is $\M_{g,1}$. Show that the isotropy group of any 
point of $\X_{g,n}$ lying above $[C;x]\in \M_{g,n}$ is naturally isomorphic
to $\Aut(C,x)$.
\end{xca}

We shall regard $\M_{g,n}$ as the orbifold $\G_{g,n}\bs \X_{g,n}$. There
is a natural isomorphism
$$
\piorb_1(\M_{g,n},[C;x_1,\dots,x_n]) \cong \pi_0 \Diff^+(C,\{x_1,\dots,x_n\})
$$
provided $\Aut(C;x_1,\dots,x_n)$ is trivial.

\section{Level Structures}

Level structures are useful technical devices for rigidifying curves.
Suppose that $\ell$ is a positive integer. A {\it level $\ell$ structure}
on a compact Riemann surface is the choice of a symplectic basis of
$H_1(C,\Z/\ell\Z)$.

\begin{xca}
Show that if $C$ is a compact Riemann surface of genus $g$, then there
is a canonical isomorphism between $H_1(C,\Z/\ell\Z)$ and the $\ell$-%
torsion points in $\Pic^0 C$. (Remark: the intersection pairing on
$H_1(C,\Z/\ell\Z)$ corresponds to the Weil pairing on the $\ell$-torsion
points.)
\end{xca}

Denote the moduli space of $n$-pointed, smooth, genus $g$ curves with a 
level $\ell$ structure by $\M_{g,n}[l]$. This can be described as a quotient
of Teichm\"uller space by a subgroup of $\G_{g,n}$ that we now describe.

Fix an $n$-pointed compact oriented reference surface $S$ of genus $g$. The
mapping class group $\G_{g,n}$ acts naturally on $H_1(S,\Z)$. Since it
preserves the intersection pairing, this leads to a homomorphism
$$
\rho : \G_{g,n} \to \Aut(H_1(S,\Z),\text{intersection form}).
$$
Define the level $\ell$ subgroup of $\G_{g,n}$ to be the kernel of the
homomorphism
$$
\rho_\ell : \G_{g,n} \to \Aut(H_1(S,\Z/\ell\Z),\text{intersection form}).
$$
This homomorphism is surjective.

\begin{xca}
Show that $\G_{g,n}[\ell]$ is a normal subgroup of finite index in $\G_{g,n}$
and that the quotient is isomorphic to $Sp_g(\Z/\ell\Z)$.
\end{xca}

\begin{xca}
Show that there is a natural bijection between $\M_{g,n}[\ell]$ and the
quotient of $\X_{g,n}$ by $\G_{g,n}[\ell]$. Show that the quotient mapping
$\M_{g,n}[\ell] \to \M_{g,n}$ that forgets the level structure has finite
degree and is Galois with Galois group $Sp_g(\Z/\ell\Z)$.
\end{xca}

Choosing a symplectic basis of $H_1(S,\Z)$ gives an isomorphism 
$$
\Aut(H_1(S,\Z),\text{intersection form}) \cong Sp_g(\Z).
$$
We therefore have a homomorphism
$$
\rho : \G_{g,n} \to Sp_g(\Z).
$$
Define the level $\ell$ subgroup $Sp_g(\Z)[\ell]$ of $Sp_g(\Z)$ to be the
kernel of the reduction mod $\ell$ homomorphism
$$
Sp_g(\Z) \to  Sp_g(\Z/\ell\Z).
$$
This homomorphism is surjective.

\begin{xca}
\begin{enumerate}
\item Show that $\G_{g,n}[\ell]$ is the inverse image of $Sp_g(\Z)[\ell]$ under
$\rho$.
\item Show that $Sp_1(\Z)$ is isomorphic to $SL_2(\Z)$ and that $\rho$ is the
standard representation when $g=1$.
\end{enumerate}
\end{xca}

\begin{theorem}[Minkowski]
The group $Sp_g(\Z)[\ell]$ is torsion free when $\ell \ge 3$. \qed
\end{theorem}

\begin{proposition}
If $2g-2+n > 0$, $g\ge 1$ and $\ell \ge 3$, then the mapping class group
$\G_{g,n}[\ell]$ is torsion free and acts fixed point freely on $\X_{g,n}$.
\end{proposition}

\begin{proof}[Sketch of Proof]
The case $g=1$ is left as an exercise. Suppose that $g\ge 2$.
We first show that $\G_{g,n}[\ell]$ acts fixed point freely on $\X_{g,n}$.
The isotropy group of a point in $\X_{g,n}$ lying over $[C;x_1,\dots,x_n]$
is isomorphic to $\Aut(C;x_1,\dots,x_n)$. This is a finite group, and is
a subgroup of $\Aut C$. It is standard that the natural representation
$\Aut C \to \Aut H^0(C,\Omega^1_C)$ is injective (Exercise: prove this. Hint:
use Riemann-Roch). 
It follows that the natural representation
$$
\Aut C \to \Aut(H_1(C,\Z),\text{intersection pairing}) \cong Sp_g(\Z).
$$ 
is injective and that $\Aut(C;x_1,\dots,x_n) \cap \G_{g,n}[\ell]$ is
trivial. (Here we are realizing $\Aut(C;x_1,\dots,x_n)$ as a subgroup of 
$\G_{g,n}$ as an isotropy group.) It follows from Minkowski's Theorem that
if $\ell \ge 3$, then $\G_{g,n}[\ell]$ acts fixed point freely on $\X_{g,n}$.

The rest of the proof is standard topology. If $\G_{g,n}[\ell]$ has a torsion 
element, then it contains a subgroup $G$ of prime order, $p$ say. Since this
acts fixed point freely on the contractible space $\X_{g,n}$, it follows that
the quotient $G\bs \X_{g,n}$ is a model of the classifying space $B(\Z/p\Z)$
of the cyclic group of order $p$. Since the model is a manifold of real
dimension $6g-6+2n$, this implies that $H^k(\Z/p\Z,\Z/p\Z)$ vanishes when
$k > 6g-6+2n$. But this contradicts the known computation that
$H^k(\Z/p\Z,\Z/p\Z)$ is non-trivial for all $k \ge 0$. The result follows.
\end{proof}

Putting together the results of this section, we have:

\begin{corollary}
If $g \ge 1$, $n\ge 0$ and $\ell \ge 3$, then $\M_{g,n}[\ell]$ is smooth and
the mapping $\X_{g,n} \to \M_{g,n}[\ell]$ is unramified with Galois group
$\G_{g,n}$; the covering
$\M_{g,n}[\ell] \to \M_{g,n}$ is a finite (ramified) Galois covering with
Galois group $Sp_g(\Z/\ell\Z)$. \qed
\end{corollary}

With this information, we are able to prove a non-trivial result about the
mapping class group.

\begin{corollary}[McCool, Hatcher-Thurston]
For all $g$ and $n$, the mapping class group is finitely presented.
\end{corollary}

\begin{proof}
We shall use the fact that each $\M_{g,n}[\ell]$ is a quasi-projective variety.
As we have seen, this is smooth when $\ell \ge 3$ and has fundamental group
isomorphic to $\G_{g,n}[\ell]$. But a well known result Lojasiewicz
\cite{lojasiewicz} (see also \cite{hironaka}) implies
that every smooth quasi-projective variety has the homotopy type of a finite
complex. It follows that $\G_{g,n}[\ell]$ is finitely presented when
$\ell \ge 3$. But since $\G_{g,n}[\ell]$ has finite index in $\G_{g,n}$,
this implies that $\G_{g,n}$ is also finitely presented.
\end{proof}

\section{Cohomology}

One way to define the homology and cohomology of a group $G$ is to find
a topological space $X$ such that
$$
\pi_j(X,*) =
\begin{cases}
G & j=1 ; \cr
0 & j\neq 1.
\end{cases}
$$
This occurs, for example, when the universal
covering of $X$ is contractible. Such a space is called a {\it classifying
space of $G$}. Under some mild hypotheses, it is unique up to homotopy.
One then defines the cohomology of $G$ with coefficients in the $G$-module
$V$ by
$$
H^j(G,V) = H^j(X,\V)
$$
where $\V$ is the locally constant sheaf over $X$ whose fiber is $V$
and whose monodromy is given by the action of $G$ on $V$. It is well
defined. Homology is defined similarly in terms of the homology of $X$
with coefficients in $\V$.

\begin{xca}
Suppose that $V$ is a $\G_{g,n}$-module and $\V$ is the corresponding
locally constant sheaf over $\M_{g,n}[\ell]$, where $\ell \ge 3$. Show
that there is a natural isomorphism
$$
H^\dot(\M_{g,n}[\ell],\V) \cong H^\dot(\G_{g,n}[\ell],V).
$$
Show that the group $Sp_g(\Z/\ell\Z)$ acts on both sides and that
the isomorphism is $Sp_g(\Z/\ell\Z)$-equivariant.
\end{xca}

Since $\G_{g,n}$ does not act fixed point freely on $\X_{g,n}$, $\M_{g,n}$
is not a classifying space for $\G_{g,n}$. Nonetheless, standard topological
arguments give
$$
H^\dot(\M_{g,n},\Q) \cong H^\dot(\M_{g,n},\Q)^{Sp_g(\Z/\ell\Z)}
$$
and
$$
H^\dot(\G_{g,n},\Q) \cong H^\dot(\G_{g,n},\Q)^{Sp_g(\Z/\ell\Z)}.
$$
It follows that:

\begin{theorem}
There is a natural isomorphism
$$
H^\dot(\G_{g,n},\Q) \cong H^\dot(\M_{g,n},\Q). \qed
$$
\end{theorem}

If $V$ is a $\G_{g,n}$-module and $\G_{g,n}$ has fixed points in $\X_{g,n}$,
we cannot always define a local system $\V$ over $\M_{g,n}$ corresponding to
$V$. However, we can formally define
$$
H^\dot(\M_g,\V\otimes\Q) :=
H^\dot(\M_{g,n}[\ell],\V\otimes\Q)^{Sp_g(\Z/\ell\Z)}.
$$
where the superscript $Sp_g(\Z/\ell\Z)$ means that we take the
$Sp_g(\Z/\ell\Z)$
invariant part. This should be regarded as the cohomology of the orbifold
$\M_{g,n}$ with coefficients in the orbifold local system corresponding to
$V$.

\begin{xca}
Show that this definition is independent of the choice of the level
$\ell \ge 3$.

\end{xca}

Let $A$ be the locus of curves in $\M_g$ with non-trivial automorphisms.
The goal of the following exercise is to show that $A$ is an analytic
subvariety of $\M_g$, each of whose components has codimension $\ge g-2$. Set
$$
\M_g' = \M_g - A.
$$

\begin{xca}
Suppose that $X$ is a Riemann surface of genus $g\ge 2$ and that $G$ is a
finite subgroup of $\Aut X$.
\begin{enumerate}
\item Set $Y=G\bs X$. Show that if $X$ is compact, then
$$
g(X) - 1 = d(g(Y) - 1) + \sum_\O(d-|\O|)/2
$$
where $\O$ ranges over the orbits of $G$ acting on $X$ and where $d$ is the
order of $G$. 
\item Show that if $g\ge 3$, then each component of the locus of curves in
$\M_g$ that have automorphisms is a proper subvariety of $\M_g$.
\item Show that the hyperelliptic locus in $\M_g$ has codimension $g-2$.
\item Show that codimension of each component of the locus in $\M_g$ with
curves with a non-trivial automorphism is $\ge g-2$, with equality if and
only if the component is the hyperelliptic locus. Hint: reduce to the
case where $d$ is prime.
\item Give an argument that there are only finitely many components of
the locus in $\M_g$ of curves with a non-trivial automorphism.
\end{enumerate}
Deduce that when $g\ge 3$, the set of points of $\M_g$ corresponding to
curves without automorphisms is Zariski dense in $\M_g$.
\end{xca}

Since each component of $A$ has real codimension $\ge 2g-4$, it follows
from standard topological arguments (transversality) that the inverse
image $\X_g'$ of $\M_g$ in $\X_g$ has the property
$$
\pi_j(\X_g') = 0 \text{ if } j<2g-5.
$$
From this, one can use standard topology to show that if $g \ge 3$, and
if $V$ is any $\G_g$-module, then there is a natural mapping
$$
H^k(\G_g,V) \to H^k(\M_g',\V)
$$
which is an isomorphism when $k < 2g-5$ and injective when $k=2g-5$.%
\footnote{There is a similar result for homology with the arrows reversed and
injectivity replaced by surjectivity.}
In particular, there is a natural isomorphism
$$
H^k(\G_g,\Z) \cong H^k(\M_g',\Z)
$$
when $k \le 2g-5$.

There are similar results for $\M_{g,n}$, but the codimension of the locus
of curves with automorphisms rises with $n$. For example, if $n > 2g+2$,
then $\Aut(C;x_1,\dots,x_n)$ is always trivial by the Lefschetz fixed point
formula. (Exercise: prove this.)

\begin{center}
{\sc Lecture 3: The Picard Group}
\end{center}

In this lecture, we compute the orbifold Picard group of $\M_g$ for all
$g\ge 1$. Recall that an orbifold line bundle over $\M_g$ is a holomorphic line
bundle $\cL$ over Teichm\"uller space $\X_g$ together with an action of
the mapping class group $\G_g$ on it such that the projection $\cL \to \X_g$
is $\G_g$-equivariant. An orbifold section of this line bundle is a holomorphic
$\G_g$-equivariant section $\X_g \to \cL$ of $\cL$. This is easily seen
to be equivalent to
fixing a level $\ell \ge 3$ and considering holomorphic line bundles over
$\M_g[\ell]$ with an $Sp_g(\Z/\ell\Z)$-action such that the projection is
$Sp_g(\Z/\ell\Z)$-equivariant. Working on $\M_g[\ell]$ has the advantage
that we can talk about algebraic line bundles more easily.

An {\it algebraic orbifold line bundle} over $\M_g$ is an algebraic
line bundle over $\M_g[\ell]$ for some $\ell$ equipped with an action
of $Sp_g(\Z/\ell\Z)$ such that the projection to $\M_g[\ell]$ is
$Sp_g(\Z/\ell\Z)$equivariant.
A section of such a line bundle is simply an $Sp_g(\Z/\ell\Z)$-equivariant
section defined over $\M_g[\ell]$. Isomorphism of such orbifold line bundles
is defined in the obvious way. Let
$$
\Pic_\orb \M_g
$$
denote the group of isomorphisms classes of algebraic orbifold line bundles
over $\M_g$. Our goal in this lecture is to compute this group. It is 
first useful to review some facts about the Picard group of a smooth
projective variety.

\section{General Facts}
Recall that if $X$ is a compact K\"ahler manifold (such as a smooth projective
variety), then the exponential sequence gives an exact sequence
$$
0 \to H^1(X,\Z) \to H^1(X,\O_X) \to \Pic X \to H^2(X,\Z)
$$
where the last map is the first Chern class $c_1$. The quotient
$$
\Pic^0 X := H^1(X,\O_X)/H^1(X,\Z)
$$
is the group of topologically trivial complex line bundles and is a compact
complex torus (in fact, an abelian variety). Note that if $H^1(X,\Z)$ vanishes,
then $\Pic^0 X$ vanishes
and the first Chern class $c_1 : \Pic X \to H^2(X,\Z)$ is injective.

If $H^1(X,\Z)$ vanishes, it follows from the Universal Coefficient Theorem
that the torsion subgroup of $H^2(X,\Z)$ is $\Hom(H_1(X,\Z),\C^\ast)$.
Every torsion element of $H^2(X,\Z)$ is the Chern class of a holomorphic
line bundle as a homomorphism $\chi : H_1(X,\Z) \to \C^\ast$ gives rise to a
flat (and therefore holomorphic) line bundle over $X$.

\begin{xca}
Suppose that $X$ is a compact K\"ahler manifold with $H^1(X,\Z) = 0$.
Show that if $\cL$ is the flat line bundle over whose monodromy is given by
the homomorphism $\chi : H_1(X,\Z) \to \C^\ast$, then
$$
c_1(\cL) = \chi \in \Hom(H_1(X,\Z),\C^\ast) \subseteq H^2(X,\Z).
$$
\end{xca}
 
These basic facts generalize to non-compact varieties. Suppose that $X$
is a smooth quasi-projective variety. Define $\Pic X$ to be the group of
isomorphism classes of algebraic line bundles over $X$, and $\Pic^0 X$ to
be the kernel of the Chern class mapping
$$
c_1 : \Pic X \to H^2(X,\Z).
$$

\begin{theorem}
Suppose that $X$ is a smooth quasi-projective variety. If $H^1(X,\Q)$ vanishes,
then $\Pic^0 X = 0$ and the torsion subgroup of $\Pic X$ is naturally
isomorphic to $\Hom(H_1(X,\Z),\C^\ast)$.
\end{theorem}

\begin{proof}[Sketch of Proof]
There are several ways to prove this. One is to use Deligne cohomology
which gives a Hodge theoretic computation of $\Pic X$. Details can
be found in \cite{hain:normal}, for example. A more elementary approach
goes as follows. First pick a smooth compactification $\Xbar$ of $X$. Each
line bundle $\cL$ over $X$ can be extended to a line bundle over $\Xbar$.
Any two extensions differ by twists by the divisors in $X$ that lie
in $\Xbar - X$. After twisting by suitable boundary components, we may
assume that the extended line bundle also has trivial $c_1$ in
$H^2(\Xbar,\Z)$. (Prove this using the Gysin sequence.)
It therefore gives an element of $\Pic^0 \Xbar$, which, by the discussion
at the beginning of the section, is flat. This implies that the
original line bundle $\cL$ over $X$ is also flat. But if $H^1(X,\Q)$
vanishes, then $H_1(X,\Z)$ is torsion and $c_1(\cL)\in H^2(X,\Z)$ is the
corresponding character $\chi : H_1(X,\Z) \to \C^\ast$. But since $c_1(\cL)$
is trivial, this implies that $\chi$ is the trivial character and that
$\cL$ is trivial.
\end{proof}

This result extends to the orbifold case. By a smooth quasi-projective
orbifold, we mean an orbifold $\G\bs X$ where $\G$ has a subgroup $\G'$
of finite index that acts fixed point freely on $X$ and where $\G'\bs X$ 
is a smooth quasi-projective variety. There is a Chern class mapping
$$
c_1 : \Pic_\orb (\G\bs X) \to H^2(\G,\Z).
$$
(The Picard group is constructed using equivariant algebraic line bundles
on finite covers of $\G\bs X$.)
The Chern class $c_1$ can be constructed using the Borel construction, for
example. Define $\Pic_\orb^0 (\G\bs X)$ to be the kernel of $c_1$.

\begin{theorem}
Suppose that $\G\bs X$ is a smooth quasi-projective orbifold. If $H^1(X,\Q)$
vanishes, then $\Pic_\orb^0 X = 0$ and the torsion subgroup of $\Pic X$ is
naturally isomorphic to $\Hom(H_1(\G,\Z),\C^\ast)$.
\end{theorem}

\begin{proof}[Sketch of Proof]
Fix a finite orbifold covering of $\G'\bs X$ of $\G\bs X$ where $\G'$ is
normal in $\G$ and acts fixed point freely on $X$. The Galois group of the
covering is $G=\G/\G'$. Suppose that $\cL$ is an orbifold line bundle over
$\G\bs X$ with $c_1(\cL) = 0$ in $H^2(\G,\Z)$. This implies that the first
Chern class of the pullback of $\cL$ to $\G'\bs X$ is also trivial. Since
this pullback has a natural $G$-action, this means that the corresponding
point in $\Pic^0 (\G'\bs X)$ is $G$-invariant. Since
$$
H^1(\G'\bs X,\Q)^G \cong H^1(\G\bs X,\Q) = 0
$$
it follows that the $G$-invariant part of $\Pic (\G'\bs X)$ is finite. It 
follows that some power $\cL^{\otimes N}$ of the pullback of $\cL$ to
$\G'\bs X$ is trivial. This also has a $G$-action. If $s$ is a
trivializing section of $ \cL^{\otimes N}$, then the product
$$
\bigotimes_{g\in G}\, g\cdot s
$$
is a $G$-invariant section of $\cL^{\otimes N|G|}$. It follows
that $\cL$ has a flat structure invariant under the $G$-action. But since
$\cL$ has trivial $c_1$, it must have trivial monodromy. It is therefore
trivial. It follows that $\Pic_\orb (\G\bs X)$ is trivial.
\end{proof}

Assembling the pieces, we have:

\begin{corollary}
\label{c_inj}
If $H^1(\G_g,\Q)$ vanishes, then the Chern class mapping
$$
\Pic_\orb \M_g \to H^2(\G_g,\Z)
$$
is injective. \qed
\end{corollary}

\section{Relations in $\G_g$}
\label{relations}

In this section we write down some well known relations that hold in
various mapping class groups. These will be enough to compute $H_1(\G_g)$,
which we shall do in the next section.

First some notation. We shall let $S$ denote any compact oriented surface
with (possibly empty) boundary. The corresponding mapping class group
is defined to be
$$
\G_S := \pi_0 \Diff^+(S,\partial S).
$$
That is, $\G_S$ consists of isotopy classes of orientation preserving 
diffeomorphisms that act trivially on the boundary $\partial S$ of $S$.

\begin{xca}
Show that elements of $\G_S$ can be represented by diffeomorphisms that equal
the identity in a neighbourhood of $\partial S$.
\end{xca}

\begin{xca}
Show that if $S$ is a compact oriented surface and $T$ is a compact
subsurface, then there is a natural homomorphism $\G_T \to \G_S$ obtained
by extending elements of $\G_T$ to be the identity outside $T$.
\end{xca}

An important special case is where we take $T$ to be the cylinder
$[0,1]\times S^1$ (with the product orientation). One element of $\G_T$ is
the diffeomorphism
$$
\tau : (t,\theta) \mapsto (t,\theta + 2\pi t).
$$
\begin{figure}[!ht]
\epsfig{file=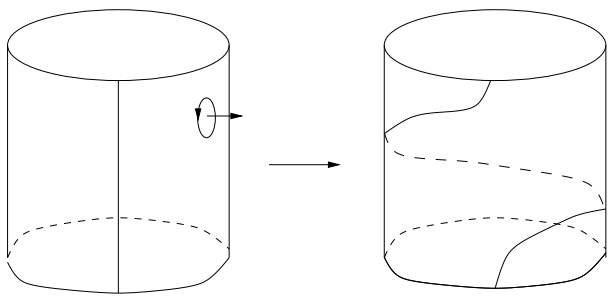, height=1.5in}
\caption{a positive Dehn twist}
\end{figure}

\begin{theorem}
If $T$ is the cylinder, then $\G_T$ is infinite cyclic and is generated
by $\tau$. \qed
\end{theorem}

\begin{figure}[!ht]
\epsfig{file=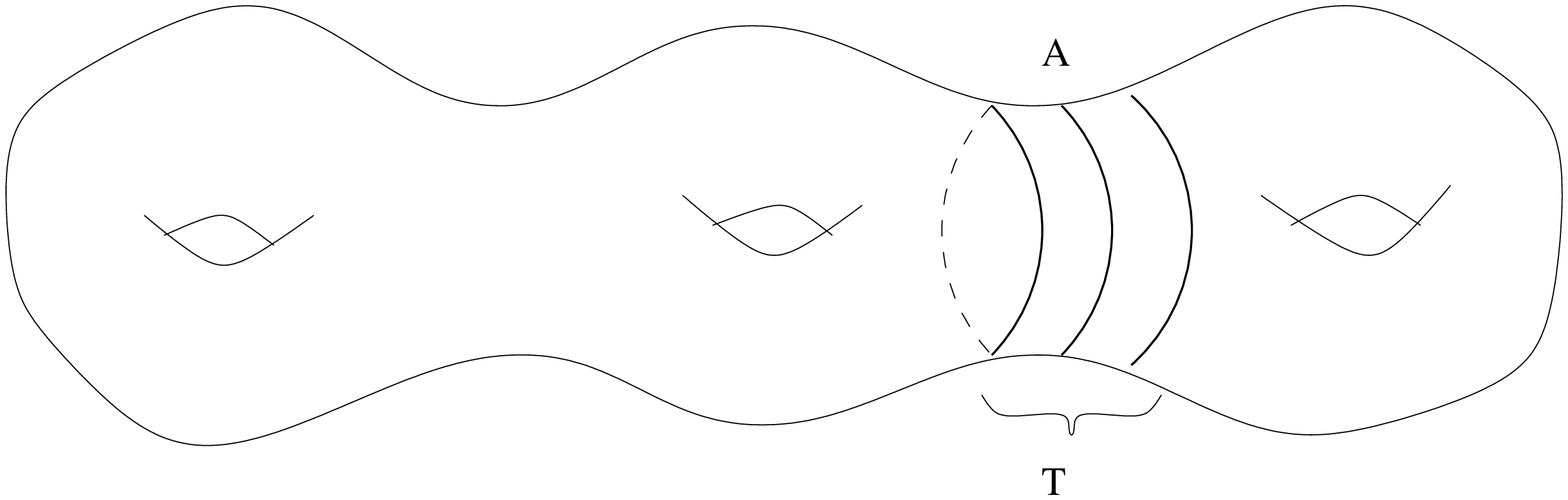, height=1in}
\caption{}
\label{collar}
\end{figure}

Now, if $S$ is any surface, and $A$ is a (smoothly imbedded) simple closed
curve (SCC), then $A$ has a neighbourhood that is diffeomorphic to the
cylinder $T = [0,1]\times S^1$ as in Figure~\ref{collar}.
Denote the image of $\tau$ under the homomorphism $\G_T \to \G_S$ by $t_A$.
It is called the {\it Dehn twist about $A$}.
\begin{theorem}
The class of $t_A$ in $\G_S$ is independent of the choice of the tubular
neighbourhood $T$ of $A$ and depends only on the isotopy class of $A$. 
Every mapping class group is generated by Dehn twists. \qed
\end{theorem}

\begin{xca}
\label{conj}
Show that if $\phi \in \G_S$ and $A$ is a SCC in $S$, then
$$
t_{\phi(A)} = \phi t_A \phi^{-1}.
$$
\end{xca}

Observe that if two elements of $\G_S$ can be represented by diffeomorphisms
with disjoint support, then they commute. In particular:

\begin{proposition}
If $A$ and $B$ are disjoint SCCs on S, then $t_A$ and $t_B$ commute in $\G_S$.
\qed
\end{proposition}

\begin{xca}
Show that if $S$ is a surface with boundary and $A$ and $B$ are each
non-separating SCCs in $S$, then there is an orientation preserving 
diffeomorphism $\phi$ of $S$ that takes one onto the other. Deduce that
$t_A$ and $t_B$ are conjugate in $\G_S$. Give an example to show that not all
Dehn twists about bounding SCCs are conjugate in $\G_S$ when $g\ge 3$.
\end{xca}

\begin{figure}[!ht]
\epsfig{file=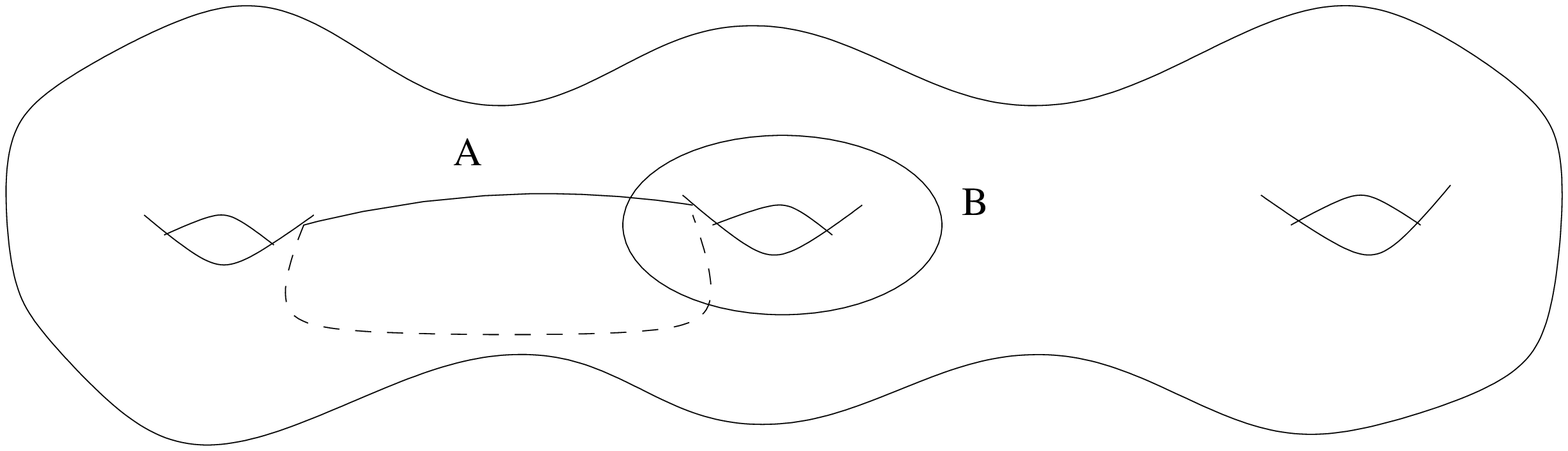, height=1in}
\caption{}
\label{scc}
\end{figure}

Next we consider what happens when $A$ and $B$ are two simple closed curves
in $S$ that meet transversally in one point as in Figure~\ref{scc}.

\begin{xca}
Show that if $A$ and $B$ are two SCCs in $S$ that meet transversally in one
point, then there is a neighbourhood of their union that is a compact genus 1
surface with one boundary component. Hint: Compute the homology of a small
regular neighbourhood of the union and then apply the classification of
compact oriented surfaces.
\end{xca}

So any relation that holds between Dehn twists such curves in a genus
one surface with one boundary component will hold in all surfaces.
Let $S$ be a compact genus 1 surface with one boundary component and
let $A$ and $B$ be the two SCCs in the diagram:
\begin{figure}[!ht]
\epsfig{file=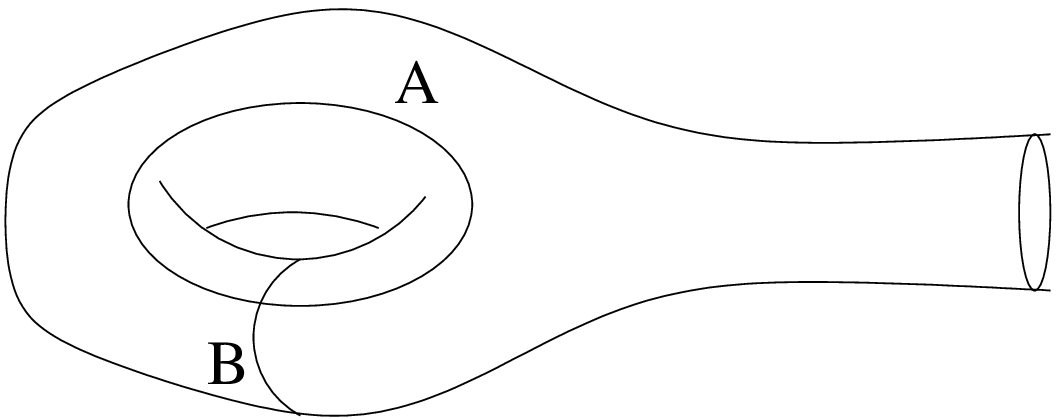, height=1in}
\caption{}
\end{figure}

\begin{theorem}
\label{braid}
With notation as above, we have
$$
t_A t_B t_A = t_B t_A t_B
$$
in $\G_S$ and therefore in all mapping class groups.
\end{theorem}

This relation is called the {\it braid relation} as it comes from the braid
group on 3 strings\footnote{A good reference for braid groups is Joan
Birman's book \cite{birman}.} using the following technique.

Denote the unit disk by $D$. We view this as a manifold with boundary.
Suppose that $P$ is a set of $n$ distinct points of $D$, none of which
lies on $\partial D$. The braid group $B_n$ is defined to be
$$
B_n := \pi_0 \Diff^+(D,(P)),
$$
where $\Diff^+(D,(P))$ denotes the group of orientation preserving 
diffeomorphisms of $D$ that fix $P$ as a set, but may permute its elements.
There is a surjective homomorphism
$$
B_n \to \Aut P \cong \Sigma_n
$$
onto the symmetric group on $n$ letters.

Suppose that $U$ is a disk imbedded in $D$ such that
$$
\partial U \cap P = \emptyset
$$
and $U \cap P$ consists of two distinct points $x$ and $y$ of
$P$. Then there is an element $\sigma_U$ of $B_n$ whose square is the Dehn
twist $t_{\partial U}$ about the boundary of $U$ and which swaps $x$ and
$y$.

\begin{figure}[!ht]
\epsfig{file=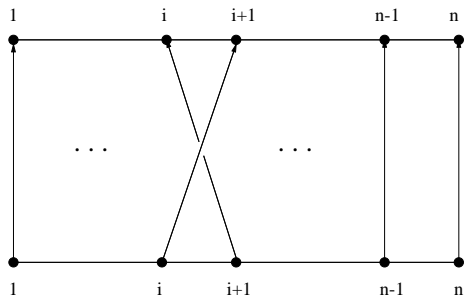, height=1.5in}
\caption{the generator $\sigma_i$}
\label{gens}
\end{figure}

It can be represented schematically as in Figure~\ref{basic_braid}:
\begin{figure}[!ht]
\epsfig{file=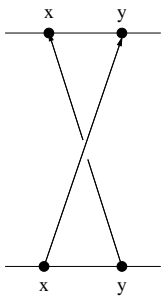, height=1.5in}
\caption{a basic braid}
\label{basic_braid}
\end{figure}
The braid group $B_n$ is generated by the braids $\sigma_1, \dots, \sigma_{n-1}$
illustrated in Figure~\ref{gens}.
Note that $\sigma_i$ and $\sigma_j$ commute when $|i-j| > 1$, and that
\begin{equation}
\label{braid_relation}
\sigma_i \sigma_{i+1}\sigma_i = \sigma_{i+1}\sigma_i\sigma_{i+1}.
\end{equation}

Now suppose that $(S,\partial S) \to (D,\partial D)$ is a branched covering,
unramified over the boundary. Suppose that the image of the branch points
is $P$. Then there is a natural homomorphism
$$
B_n \to \G_S.
$$
So relations in $B_n$ will give relations in $\G_S$. The relations we
are interested in come from certain double branched coverings of the disk.

\begin{figure}[!ht]
\epsfig{file=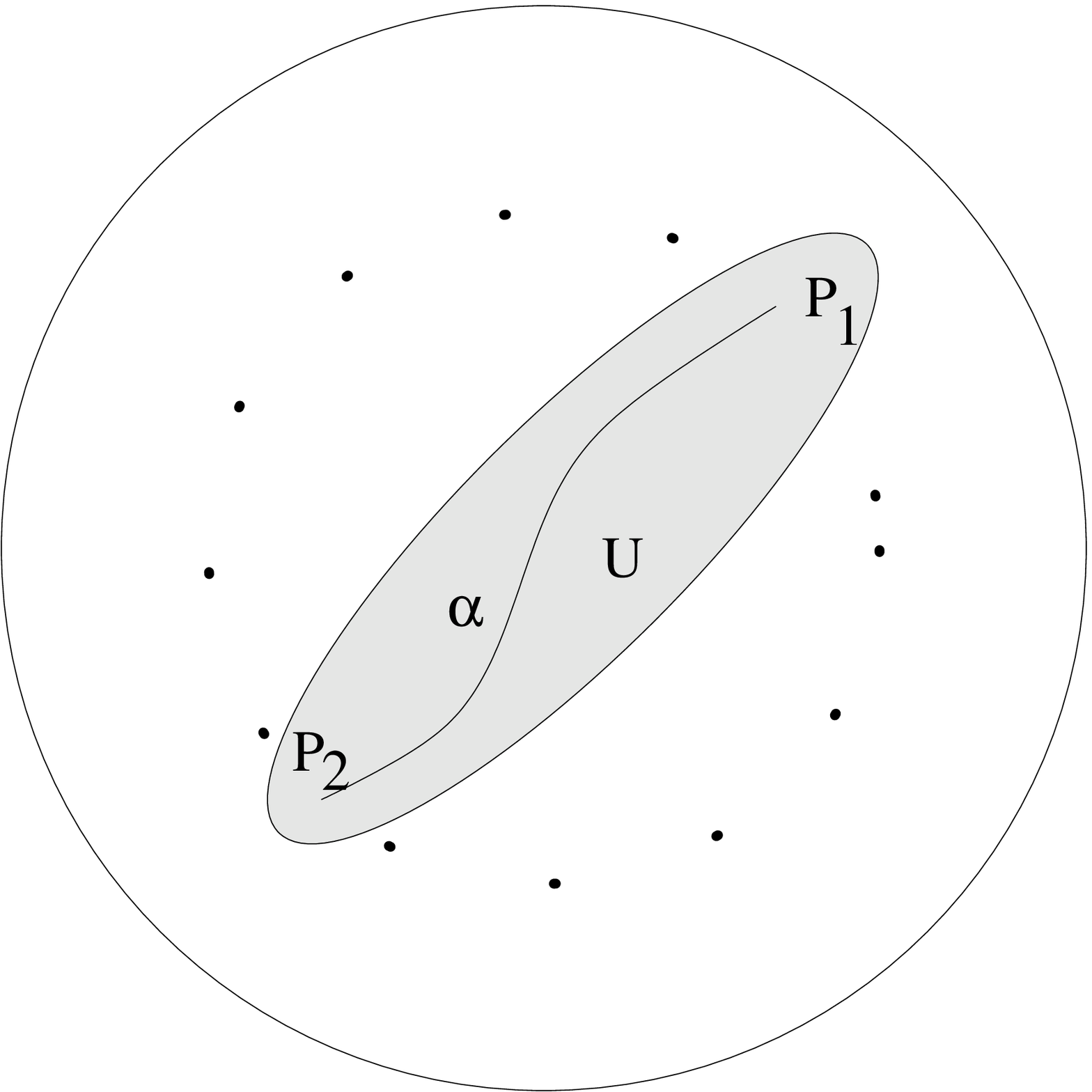, height=2in}
\caption{}
\label{braid_pic}
\end{figure}

Suppose that $(S,\partial S) \to (D, \partial D)$ is a 2-fold branched covering.
The inverse image of a smooth arc $\alpha$ joining two critical values
$p_1,p_2 \in P$ and avoiding $P$ otherwise, is an SCC, say A, in $S$. There
is  a small neighbourhood $U$ of $\alpha$ that is diffeomorphic to a disk and
whose intersection with $P$ is $\{p_1,p_2\}$ as in Figure~\ref{braid_pic}.
There is a braid $\sigma$ supported in $U$ that swaps $p_1$ and $p_2$ by
rotating in the positive direction about $\alpha$.

\begin{proposition}
The image of $\sigma$ under the homomorphism $B_n \to \G_S$ is the Dehn
twist $t_A$ about $A$. \qed
\end{proposition}

\begin{xca}
Show that a genus 1 surface $S$ with one boundary component can be realized
as a 2:1 covering of $D$, branched over 3 points.
\end{xca}

\begin{figure}[!ht]
\epsfig{file=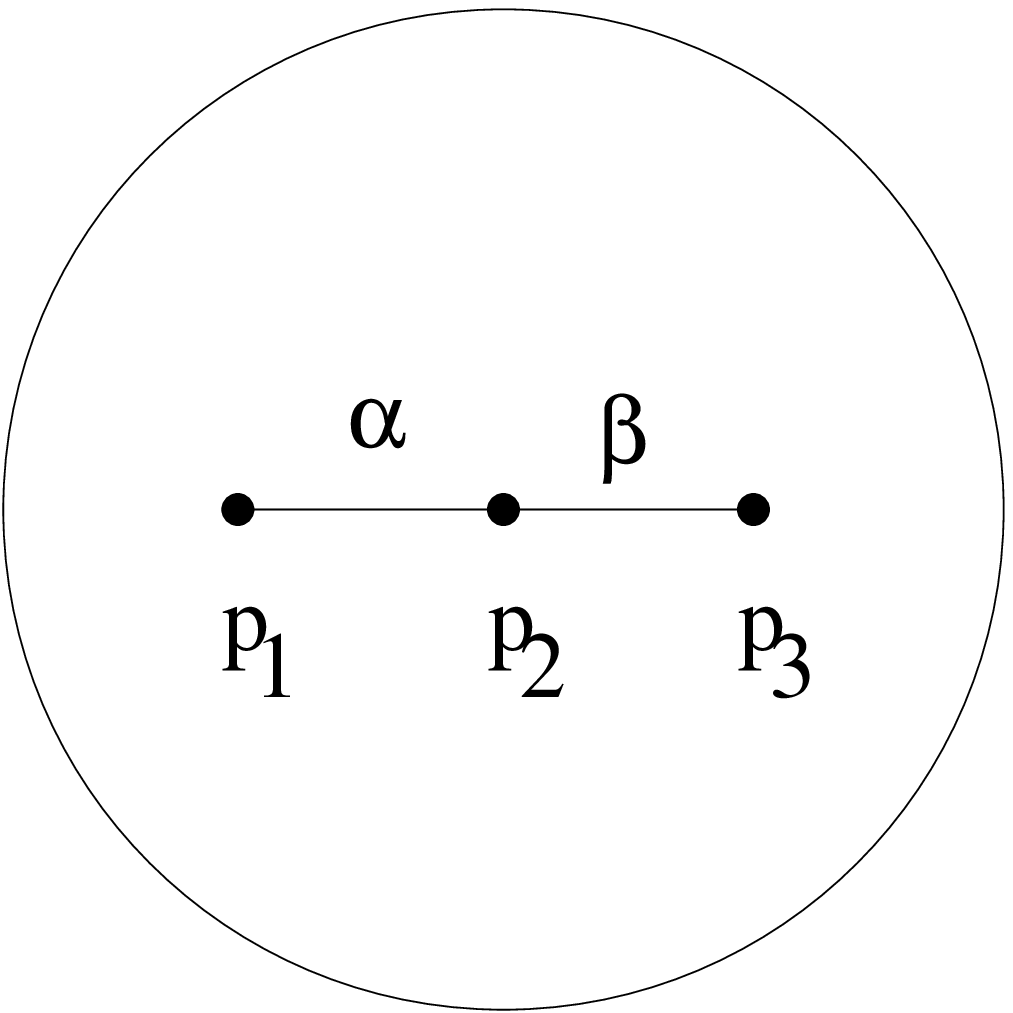, height=1.5in}
\caption{}
\label{slits}
\end{figure}

We therefore have a homomorphism $B_3 \to \G_S$. Note that the
inverse image of the two arcs $\alpha$ and $\beta$ in Figure~\ref{slits}
under the covering of the disk branched over $\{p_1,p_2,p_3\}$ is a pair of
SCCs in $S$ that intersect transversally in one point. The braid relation in
$\G_S$ now follows as we have the braid relation (\ref{braid_relation})
$$
\sigma_1 \sigma_2 \sigma_1 = \sigma_2 \sigma_1 \sigma_2
$$
in $B_3$.

More relations can be obtained this way. Suppose that $S$ is a compact
genus 1 surface with 2 boundary components.
\begin{figure}[!ht]
\epsfig{file=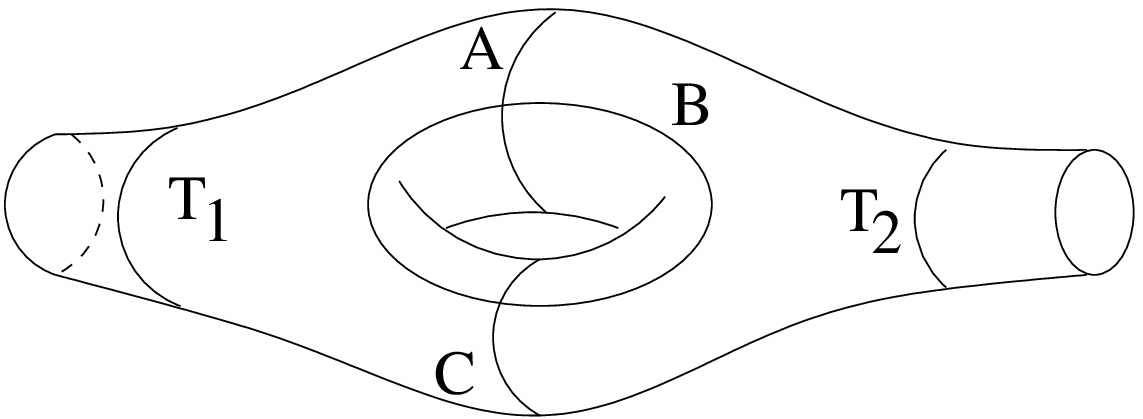, height=1in}
\caption{}
\label{two_holes}
\end{figure}

Let $A$, $B$, $C$, $T_1$ and $T_2$ be the SCCs in Figure~\ref{two_holes}.
Denote the Dehn twists about $A$, $B$ and $C$ by $a$, $b$ and $c$, and those
about $T_1$ and $T_2$ by $t_1$ and $t_2$.

\begin{theorem}
\label{2hole}
The relation
$$
(abc)^4 = t_1 t_2
$$
holds in the mapping class group $\G_S$.
\end{theorem}

\begin{corollary}
\label{1hole}
If $S$ is a compact genus 1 surface with one boundary component, then
the relation
$$
(ab)^6 = t
$$
holds in $\G_S$, where $a$ and $b$ denote Dehn twists about a pair of SCCs
that intersect transversally in one point and $t$ denotes a Dehn twist
about the boundary.
\end{corollary}

\begin{xca}
Deduce Corollary~\ref{1hole} from Theorem~\ref{2hole} by capping off one
boundary component and using the braid relation, Theorem~\ref{braid}. 
\end{xca}

Theorem~\ref{2hole} is proved in the following exercise.

\begin{xca}
Suppose that $S$ is a genus 1 surface with two boundary components.
Show that $S$ is a double covering of the disk, branched over 4 points. Use
this to construct a homomorphism from the braid group $B_4$
on 4 strings into $\G_S$.
Use this to prove the relation
$$
(abc)^4 = t_1 t_2
$$
where $a$, $b$, $c$, $t_1$ and $t_2$ denote Dehn twists on the SCCs $A$, $B$,
$C$, $T_1$ and $T_2$ in the diagram above.
Note that in the braid group, we have the relation
$$
t = \sigma_1\sigma_2\sigma_3
$$
where $t$ is Dehn twist about the boundary of the disk and $\sigma_i$ is the
$i$th standard generator of $B_4$.
\end{xca}

\begin{figure}[!ht]
\epsfig{file=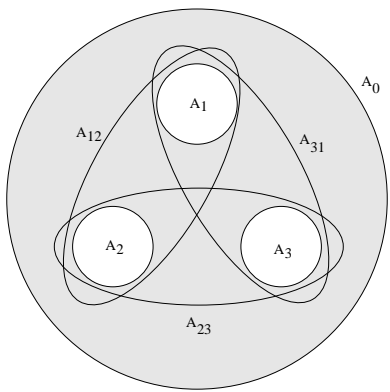, height=2in}
\caption{the lantern configuration}
\label{lantern}
\end{figure}

There is one final relation. It is called the {\it lantern relation}
and is due to Johnson and Harer independently. Let $S$ be a disk with
3 holes. (That is, a genus 0 surface with 4 boundary components.)
Consider the SCCs in Figure~\ref{lantern}.

\begin{theorem}
The relation
$$
a_0 a_1 a_2 a_3 = a_{12} a_{23} a_{31}
$$
holds in $\G_S$ where $a_i$ denotes the Dehn twist about $A_i$ and
$a_{ij}$ denotes the Dehn twist about $A_{ij}$. \qed
\end{theorem}

\section{Computation of $H_1(\G_g,\Z)$}

The fact $\G_g$ is generated by Dehn twists and the relations given in
the previous section allow the computation of $H_1(\G_g)$.

\begin{theorem}[Harer]
If $g\ge 1$, then
$$
H_1(\G_g,\Z) \cong
\begin{cases}
\Z/12\Z & g=1; \cr
\Z/10\Z & g=2; \cr
0 	& g \ge 3.
\end{cases}
$$
\end{theorem}

\begin{proof}
We begin with the observation that if $S$ is a compact oriented genus $g$
surface, then all Dehn twists on non-separating SCCs lie in the same homology
class as they are conjugate by Exercise~\ref{conj}. Denote their common
homology class by $L$. Next, using the relations coming from an imbedded
genus 1 surface with one boundary component, we see that the homology
class of any separating SCC that divides $S$ into a genus 1 and genus $g-1$
surfaces has homology class $12L$. Using the relation that comes from
an imbedded genus 1 surface with 2 boundary components, we see that
the homology class of every separating SCC is an integer multiple of $L$. It
follows that $H_1(\G_g)$ is cyclic and generated by $L$.

Now suppose that $g\ge 3$. Then we can find an imbedded lantern as
in Figure~\ref{imbedded}
\begin{figure}[!ht]
\epsfig{file=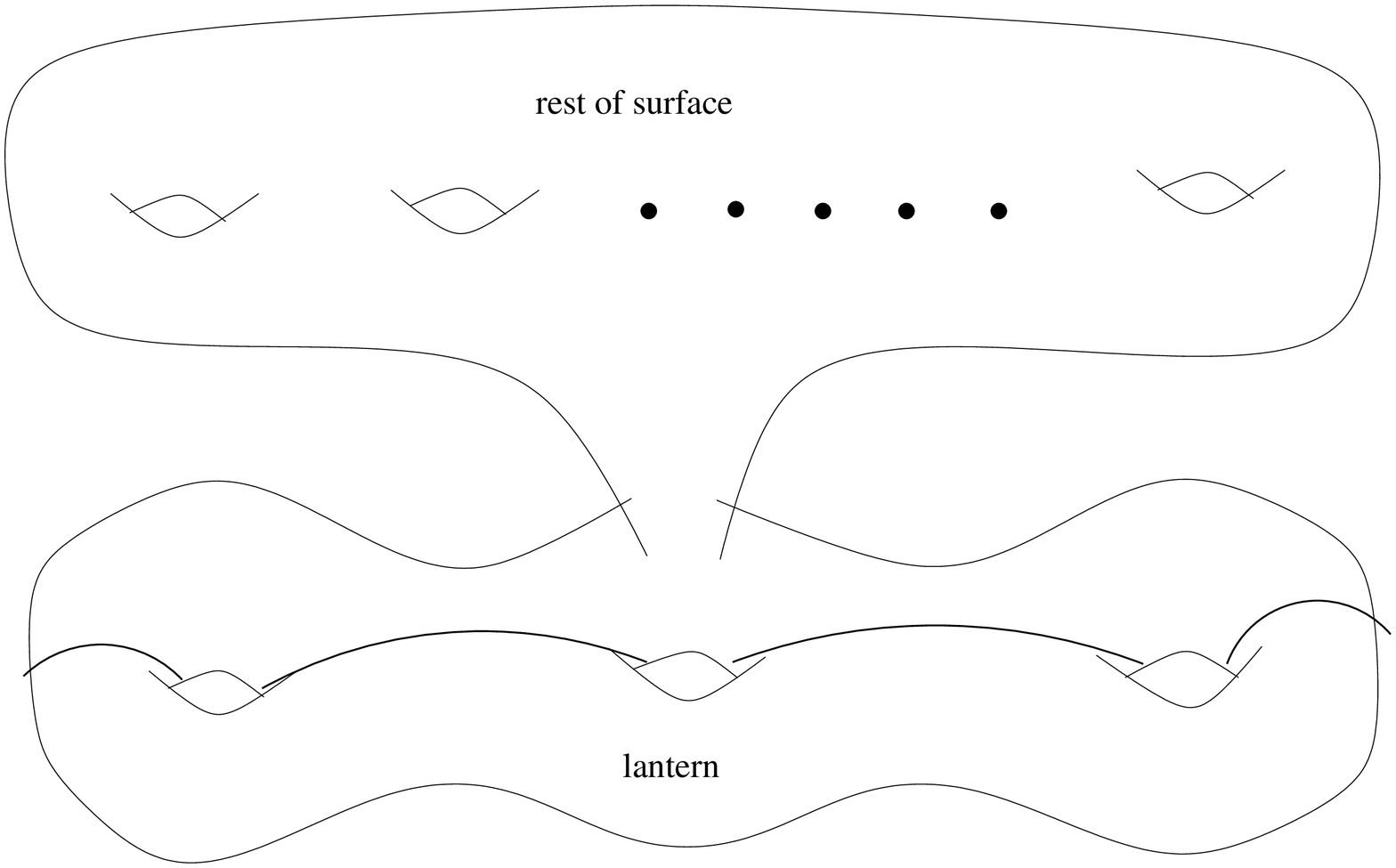, height=2.5in}
\caption{}
\label{imbedded}
\end{figure}
Since each of the curves in this lantern is non-separating, the lantern
relation tells us that $3L = 4L$, which implies that $L=0$. This proves
the vanishing of $H_1(\G_g)$ when $g\ge 3$.

When $g=1$, the relation for a genus 1 surface with one boundary component
implies that $12L = 0$ as the twist about the boundary is trivial in $\G_1$.
Thus $H_1(\G_1)$ is a quotient of $\Z/12\Z$. But, as we shall explain a little
later, the fact that $\Pic_\orb \M_1$ is at least as big as $\Z/12$ implies
that $H_1(\G_1) = \Z/12\Z$.\footnote{The abelianization of $SL_2(\Z)$ can also
be computed using, for example, the presentation of $PSL_2(\Z)$ given in
\cite{serre}.}

A genus 2 surface can be obtained from a genus 1 surface with two boundary
components by identifying the boundary components. The relation obtained for
a genus 1 surface with
2 boundary components gives $12L = 2 L$, so that $10L = 0$. This shows
that $H_1(\G_2)$ is a quotient of $\Z/10\Z$. But, as in the genus 1 case,
the theory of Siegel modular forms shows that it cannot be any smaller.
So we have $H_1(\G_2) = \Z/10\Z$.
\end{proof}

\begin{corollary}
For all $g\ge 1$, $H^1(\G_g,\Z)$ vanishes. \qed
\end{corollary}

\section{Computation of $\Pic_\orb \M_g$}

Since $H^1(\G_g,\Q)$ is torsion, it follows from Corollary~\ref{c_inj}
that
$$
c_1 : \Pic_\orb \M_g \to H^2(\G_g,\Z)
$$
is injective. Since $H_1(\G_g)$ is torsion, the Universal Coefficient Theorem
implies that the sequence
$$
0 \to \Hom(H_1(\G_g,\Z),\C^\ast) \to H^2(\G_g,\Z) \to \Hom(H_2(\G_g),\Z) \to 1
$$
is exact. To determine the rank of $H^2(\G_g)$, we need the following
fundamental and difficult result of Harer \cite{harer:h2}.

\begin{theorem}
The rank of $H_2(\G_g,\Q)$ is $0$ if $g \le 2$ and $1$ if $g\ge 3$. \qed
\end{theorem}

Combining this with our previous discussion, we have:

\begin{theorem}
If $g\ge 1$, then
$$
H^2(\G_g,\Z) \cong
\begin{cases}
\Z/12\Z & g=1; \cr
\Z/10\Z & g=2; \cr
\Z & g \ge 3.
\end{cases}
$$
\end{theorem}

There is one obvious orbifold line bundle over $\M_g$. Namely, if we take
the universal curve $\pi : \cC \to \M_g$ over the orbifold $\M_g$, then we
can form the line bundle
$$
\cL := \det \pi_\ast \omega_{\cC/\M_g}
$$
over $\M_g$, where $\omega_{\cC/\M_g}$ denotes the relative dualizing sheaf.
The fiber over $[C]\in \M_g$ is
$$
\Lambda^g H^0(C,\Omega^1_C).
$$

\begin{theorem}
If $g\ge 1$, the orbifold Picard group of $\M_g$ is cyclic, generated by
$\cL$ and given by
$$
\Pic_\orb \M_g \cong
\begin{cases}
\Z/12\Z & g=1; \cr
\Z/10\Z & g=2; \cr
\Z & g \ge 3.
\end{cases}
$$
\end{theorem}

\begin{proof}
All but the generation by $\cL$ follows from preceding results. In genus
1 and 2, the generation by $\cL$ follows from the theory of modular forms.

Suppose that $g\ge 3$. Denote the
first Chern class of $\cL$ by $\lambda$. To prove that $\cL$ generates
$\Pic_\orb \M_g$, it suffices to show that $\lambda$ generates $H^2(\G_g,\Z)$.
The following proof of this I learned from Shigeyuki Morita. It assumes
some knowledge of characteristic classes.  A good reference for this topic
is the book \cite{milnor-stasheff} by Milnor and Stasheff.

We begin by recalling the definition of the signature of a compact oriented
4-manifold. Every symmetric bilinear form on a real vector space can be
represented by a symmetric matrix. The signature of a non-degenerate
symmetric bilinear form is the number of positive eigenvalues of a representing
matrix minus the number of negative eigenvalues. To each compact
oriented 4-manifold $X$, we associate the symmetric bilinear form
$$
H^2(X,\R) \otimes H^2(X,\R) \to \R
$$
defined by
$$
\xi_1 \otimes \xi_2 \mapsto \int_X \xi_1 \wedge \xi_2.
$$
Poincar\'e duality implies that it is non-degenerate. The signature
$\tau(X)$ of $X$ is defined to be the signature of this bilinear form. It is
a cobordism invariant.

The Hirzebruch Signature Theorem (see \cite[Theorem~19.4]{milnor-stasheff},
for example) asserts that
$$
\tau(X) = \frac{1}{3} \int_X p_1(X)
$$
where $p_1(X) \in H^4(X,\Z)$ is the first Pontrjagin class of $X$. When
$X$ is a complex manifold, 
\begin{equation}
\label{pont}
p_1(X) = c_1(X)^2 - 2 c_2(X).
\end{equation}

Now suppose that $X$ is a smooth algebraic surface and that $T$ is a smooth
algebraic curve. Suppose that $\pi : X \to T$ is a family whose fibers
are smooth curves of genus $g \ge 3$.\footnote{Such families exist --- see
\cite[p.~45, p.~55]{harris-morrison}.} Denote the relative cotangent bundle
$\w_{X/T}$ of $\pi$ by $\w$. Then it follows from the exact
sequence
$$
0 \to \pi^\ast \Omega^1_S \to \Omega^2_X \to \w_{X/T} \to 0
$$
that
$$
c_1(X) = \pi^\ast c_1(T) - c_1(\w) \text{ and }
c_2(X) = - c_1(\w)\wedge \pi^\ast c_1(T).
$$
Plugging these into (\ref{pont}) we see that $p_1(X) = c_1(\w)^2$.
Using integration over the fiber, we have
$$
\tau(X) = \frac{1}{3}\int_X c_1(\w)^2 = \frac{1}{3}\int_T \pi_\ast (c_1(\w)^2).
$$
It is standard to denote $\pi_\ast (c_1(\w)^2)$ by $\kappa_1$. Thus we have
$$
\tau(X) = \frac{1}{3}\int_T \kappa_1.
$$

An easy consequence of the Grothendieck-Riemann-Roch Theorem is that
$\kappa_1 = 12 \lambda$.
This is proved in detail in the book of Harris and Morrison
\cite[pp.~155--156]{harris-morrison}. It follows that for a family
$\pi : X \to T$ of smooth curves
$$
\tau(X) = 4 \int_T \lambda.
$$

The last step is topological.
Suppose that $F$ is a compact oriented surface and that $p : W \to F$ is an
oriented surface bundle over $F$ where the fibers of $p$ are compact oriented
surfaces of genus $g \ge 3$. Denote the local system of the first integral
homology groups of the fibers by $\H$. There is a symmetric bilinear form
\begin{equation}
\label{pairing}
H^1(F,\H_\R)\otimes H^1(F,\H_\R) \to \R
\end{equation}
obtained from the cup product and the intersection form. 
Poincar\'e duality implies that it is non-degenerate. It follows from the
Leray-Serre spectral sequence of $p$ that
$$
\tau(W) = - \text{ the signature of the pairing (\ref{pairing})}.
$$

The local system $\H$ over $F$ corresponds to a mapping
$\phi : F \to BSp_g(\Z)$ into the classifying space of the symplectic group.
Meyer \cite{meyer} shows that there is a cohomology class
$m \in H^2(Sp_g(\Z),\Z)$ whose value on $\phi$ is the signature of the
pairing (\ref{pairing}). It follows from this and the discussion above,
that under the mapping
$$
\rho^\ast : H^2(Sp_g(\Z),\Z) \to H^2(\G_g,\Z)
$$
induced by the canonical homomorphism $\rho : \G_g \to Sp_g(\Z)$, $m$ goes
to $-4\lambda$. Mayer also shows that the image of
$$
m : H^2(\G_g,\Z) \to \Z
$$
is exactly $4\Z$, which implies that $\lambda$ generates
$H^2(\G_g,\Z)$.
\end{proof}

As a corollary of the proof, we have:

\begin{corollary}
For all $g\ge 3$, both $H^2(Sp_g(\Z),\Z)$ and $H^2(\G_g,\Z)$ are
generated by $\lambda$ and the natural mapping
$$
\rho^\ast : H^2(Sp_g(\Z),\Z) \to H^2(\G_g,\Z)
$$
is an isomorphism. \qed
\end{corollary}

Note that $H^2(Sp_2(\Z),\Z)$ is infinite cyclic, while
$H^2(\G_2,\Z)$ is cyclic of order 10. So $\rho^\ast$ is not an isomorphism
in genus 2.



\end{document}